\documentstyle[11pt]{article}
\makeatletter
\ifx\@auxdefsloaded\relax\endinput
\else\let\@auxdefsloaded\relax\fi

\def\newenvironment{%
   \@ifnextchar *{\@@newenv{\global\@ignoretrue}}{\@@newenv{}*}}

\def\@@newenv#1*#2{%
   \@ifnextchar [{\@newenv{#1}{#2}}{\@newenv{#1}{#2}[0]}}

\long\def\@newenv#1#2[#3]#4#5{%
   \expandafter\newcommand\csname#2\endcsname[#3]{#4}%
   \expandafter\long\expandafter\def\csname end#2\endcsname{#5#1}}

\def\renewenvironment{%
   \@ifnextchar *{\@@renewenv{\global\@ignoretrue}}{\@@renewenv{}*}}

\def\@@renewenv#1*#2{%
   \@ifnextchar [{\@renewenv{#1}{#2}}{\@renewenv{#1}{#2}[0]}}

\long\def\@renewenv#1#2[#3]#4#5{%
   \expandafter\renewcommand\csname#2\endcsname[#3]{#4}%
   \expandafter\long\expandafter\def\csname end#2\endcsname{#5#1}}

\def\newoptcommand#1#2{%
   \@ifnextchar [{\@optargdef#1#2}{\@optargdef#1#2[1]}}

\def\renewoptcommand#1#2{%
   \edef\@tempa{\expandafter\@cdr\string#1\@nil}%
   \@ifundefined{\@tempa}{%
      \@latexerr{\string#1\space undefined}\@ehc}{}%
   \@ifnextchar [{\@reoptargdef#1#2}{\@reoptargdef#1#2[1]}}

\long\def\@optargdef#1#2[#3]#4{%
   \@ifdefinable #1{\@reoptargdef#1#2[#3]{#4}}}

\long\def\@reoptargdef#1#2[#3]#4{%
   \@tempcnta#3\relax \@tempcntb \@ne
   \let#1\relax \let\@tempa\relax
   \edef\@tempb{\long\def\csname\string#1\endcsname
      [\@tempa\the\@tempcntb]}%
   \advance\@tempcntb \@ne \advance\@tempcnta \m@ne
   \@whilenum\@tempcnta>0\do{%
      \edef\@tempb{\@tempb\@tempa\the\@tempcntb}%
      \advance\@tempcntb \@ne \advance\@tempcnta \m@ne}%
   \let\@tempa=##\@tempb{#4}%
   \def#1{\@ifnextchar [{\csname\string#1\endcsname}{%
      \csname\string#1\endcsname[#2]}}}

\def\newoptenvironment{%
   \@ifnextchar *{\@@newoptenv{\global\@ignoretrue}}{%
      \@@newoptenv{}*}}

\def\@@newoptenv#1*#2#3{%
   \@ifnextchar [{\@newoptenv{#1}{#2}{#3}}{%
      \@newoptenv{#1}{#2}{#3}[0]}}

\long\def\@newoptenv#1#2#3[#4]#5#6{%
   \expandafter\newoptcommand\csname#2\endcsname{#3}[#4]{#5}%
   \expandafter\long\expandafter\def\csname end#2\endcsname{#6#1}}

\def\renewoptenvironment{%
   \@ifnextchar *{\@@renewoptenv{\global\@ignoretrue}}{%
      \@@renewoptenv{}*}}

\def\@@renewoptenv#1*#2#3{%
   \@ifnextchar [{\@renewoptenv{#1}{#2}{#3}}{%
      \@renewoptenv{#1}{#2}{#3}[0]}}

\long\def\@renewoptenv#1#2#3[#4]#5#6{%
   \expandafter\renewoptcommand\csname#2\endcsname{#3}[#4]{#5}%
   \expandafter\long\expandafter\def\csname end#2\endcsname{#6#1}}

\newcounter{keepoptional}
\newcounter{optnestctr}

\newoptenvironment*{optional}{1}[1]{%
   \ifnum#1>\value{keepoptional}\relax
      \setcounter{optnestctr}{0}\@powerup\expandafter\@endopt\fi
   \ignorespaces}{}

\def\@powerup{\catcode`\{=12 \catcode`\}=12 \catcode`\\=12 \relax}
\def\@powerdown{\catcode`\{=1 \catcode`\}=2 \catcode`\\=0 \relax}
\begingroup \catcode`|=0 \catcode`[=1 \catcode`]=2 \@powerup
   |long|gdef|@endopt#1\end{optional}[%
      |@beginopt#1\begin{optional}*]%
   |long|gdef|@beginopt#1\begin{optional}[%
      |@ifstar[%
         |ifnum|value[optnestctr]>0|relax
            |addtocounter[optnestctr][-1]|let|@zzzap=|@endopt
         |else
            |@powerdown|end[optional]|let|@zzzap=|relax|fi
	 |@zzzap]%
      [
         |addtocounter[optnestctr][1]|@beginopt]]%
|endgroup

\makeatother
\makeatletter
\ifx\@auxdefsloaded\relax
\else \input{auxdefs.sty}\fi
\newskip\dgARROWLENGTH  \dgARROWLENGTH=2.5em\relax
\newskip\dgHORIZPAD     \dgHORIZPAD=1em\relax
\newskip\dgVERTPAD      \dgVERTPAD=2ex\relax
\newskip\dgLABELOFFSET  \dgLABELOFFSET=.7ex\relax
\newcount\dgARROWPARTS  \dgARROWPARTS=4\relax
\newcommand{\dgeverynode}{\displaystyle}
\newcommand{\dgeverylabel}{\scriptstyle}
\newskip\dgDOTSPACING   \dgDOTSPACING=0.35em
\newskip\dgDOTSIZE      \dgDOTSIZE=1.5\fontdimen8\tenln
\newcount\dgMAXSQUARE   \dgMAXSQUARE=4\relax
\newcount\dgMINDBLSQ    \dgMINDBLSQ=10\relax
\newskip\dgCOLUMNWIDTH  \dgCOLUMNWIDTH=2em\relax
\chardef\f@ur=4
\@namedef{dgo@}{\let\dg@VECTOR=\vector
   \dg@LBLPOS=\dgARROWPARTS \divide\dg@LBLPOS\tw@}%
\@namedef{dgo@-}{\let\dg@VECTOR=\line}%
\@namedef{dgo@!}{\let\dg@VECTOR=\dg@novector}%
\@namedef{dgo@1}{\dg@LBLPOS=\@ne\relax}%
\@namedef{dgo@2}{\dg@LBLPOS=\tw@\relax}%
\@namedef{dgo@3}{\dg@LBLPOS=\thr@@\relax}%
\@namedef{dgo@4}{\dg@LBLPOS=\f@ur\relax}%
\@namedef{dgo@5}{\dg@LBLPOS=5\relax}%
\@namedef{dgo@6}{\dg@LBLPOS=6\relax}%
\@namedef{dgo@7}{\dg@LBLPOS=7\relax}%
\@namedef{dgo@8}{\dg@LBLPOS=8\relax}%
\@namedef{dgo@9}{\dg@LBLPOS=9\relax}%
\def\dgt@e{\dg@DX=\@ne \dg@DY=\z@ \dg@SIZE=\@ne}%
\def\dgt@w{\dg@DX=\m@ne \dg@DY=\z@ \dg@SIZE=\@ne}%
\def\dgt@n{\dg@DX=\z@ \dg@DY=\@ne \dg@SIZE=\@ne}%
\def\dgt@s{\dg@DX=\z@ \dg@DY=\m@ne \dg@SIZE=\@ne}%
\def\dgt@ne{\dg@DX=\@ne \dg@DY=\@ne \dg@SIZE=\@ne}%
\def\dgt@se{\dg@DX=\@ne \dg@DY=\m@ne \dg@SIZE=\@ne}%
\def\dgt@nw{\dg@DX=\m@ne \dg@DY=\@ne \dg@SIZE=\@ne}%
\def\dgt@sw{\dg@DX=\m@ne \dg@DY=\m@ne \dg@SIZE=\@ne}%
\def\dgt@nne{\dg@DX=\@ne \dg@DY=\tw@ \dg@SIZE=\@ne}%
\def\dgt@nnw{\dg@DX=\m@ne \dg@DY=\tw@ \dg@SIZE=\@ne}%
\def\dgt@sse{\dg@DX=\@ne \dg@DY=-\tw@ \dg@SIZE=\@ne}%
\def\dgt@ssw{\dg@DX=\m@ne \dg@DY=-\tw@ \dg@SIZE=\@ne}%
\def\dgt@ene{\dg@DX=\tw@ \dg@DY=\@ne \dg@SIZE=\tw@}%
\def\dgt@ese{\dg@DX=\tw@ \dg@DY=\m@ne \dg@SIZE=\tw@}%
\def\dgt@wnw{\dg@DX=-\tw@ \dg@DY=\@ne \dg@SIZE=\tw@}%
\def\dgt@wsw{\dg@DX=-\tw@ \dg@DY=\m@ne \dg@SIZE=\tw@}%
\def\dggeometry{
   \dg@ZTEMP=\dg@XGRID \multiply\dg@ZTEMP\tw@
   \ifnum\dg@YGRID=\z@ \dg@ZTEMP=\tw@
   \else \divide\dg@ZTEMP\dg@YGRID \fi
   \ifnum\dg@ZTEMP>\f@ur \dg@ZTEMP=\f@ur \fi
   \ifnum\dg@ZTEMP<\@ne \dg@ZTEMP=\@ne \fi
   \unitlength=2sp\relax
   \ifnum\dg@ZTEMP<\tw@
      \advance\dg@ZTEMP\@ne
      \multiply\unitlength\dg@YGRID
   \else
      \multiply\unitlength\dg@XGRID \divide\unitlength\dg@ZTEMP
   \fi
   \dg@XGRID=\dg@ZTEMP \dg@YGRID=\tw@
   \dg@rmcommondiv\tw@\dg@XGRID\dg@YGRID
   \divide\unitlength\dg@YGRID \divide\unitlength\@m\relax}
\@namedef{dgo@..}{\let\dg@VECTOR=\dg@dotvector}%

\def\dg@dotvector(#1,#2)#3{%
   \begingroup
   \dg@XTEMP=#1\relax \dg@YTEMP=#2\relax
   \let\dg@NDOTS=\dg@XEND \let\dg@DOTDIAM=\dg@WEND
   \dg@NDOTS=\unitlength \multiply\dg@NDOTS #3\relax
   \dg@ZTEMP=\dg@YTEMP \dg@changesign\dg@YTEMP\dg@ZTEMP
   \ifnum\dg@XTEMP>\z@
      \ifnum\dg@YTEMP>\dg@XTEMP
         \multiply\dg@NDOTS\dg@YTEMP \divide\dg@NDOTS\dg@XTEMP \fi
   \else\ifnum\dg@XTEMP<\z@
      \ifnum\dg@YTEMP>-\dg@XTEMP
         \multiply\dg@NDOTS\dg@YTEMP \divide\dg@NDOTS-\dg@XTEMP \fi
   \fi\fi
   \dg@YTEMP=\dg@ZTEMP
   \divide\dg@NDOTS\dgDOTSPACING
   \ifnum\dg@NDOTS>\z@\else \dg@NDOTS=\@ne \fi
   \dg@ZTEMP=\unitlength \multiply\dg@ZTEMP #3\relax
   \divide\dg@ZTEMP\dg@NDOTS
   \ifnum\dg@XTEMP=\z@
      \dg@changesign\dg@ZTEMP\dg@YTEMP \dg@YTEMP=\dg@ZTEMP
   \else
      \dg@changesign\dg@ZTEMP\dg@XTEMP
      \multiply\dg@YTEMP\dg@ZTEMP \divide\dg@YTEMP\dg@XTEMP
      \dg@XTEMP=\dg@ZTEMP
   \fi
   \divide\dg@XTEMP\unitlength \divide\dg@YTEMP\unitlength
   \begin{picture}(0,0)
      \dg@DOTDIAM=\dgDOTSIZE \divide\dg@DOTDIAM\unitlength
      \multiput(0,0)(\dg@XTEMP,\dg@YTEMP){\dg@NDOTS}{%
         \circle*{\dg@DOTDIAM}}%
      \multiply\dg@XTEMP\dg@NDOTS \multiply\dg@YTEMP\dg@NDOTS
      \put(\dg@XTEMP,\dg@YTEMP){\vector(#1,#2){0}}%
   \end{picture}%
   \endgroup}%
\newif\ifdg@LATEXGEOM \dg@LATEXGEOMfalse

\@ifundefined{lamsvector}{}{%
   \@namedef{dgo@}{%
      \let\dg@VECTOR=\lamsvector
      \lamsshaft{?}\lamssource{?}\lamstarget{?}%
      \dg@LBLPOS=\dgARROWPARTS \divide\dg@LBLPOS\tw@\relax}%
   \@namedef{dgo@..}{\lamsshaft{.}}%
   \@namedef{dgo@=}{\lamsshaft{=}}%
   \@namedef{dgo@-}{\lamstarget{0}}%
   \@namedef{dgo@'}{\lamstarget{'}}%
   \@namedef{dgo@`}{\lamstarget{`}}%
   \@namedef{dgo@A}{\lamstarget{A}}%
   \@namedef{dgo@V}{\lamssource{V}}%
   \@namedef{dgo@J}{\lamssource{J}}%
   \@namedef{dgo@L}{\lamssource{L}}%
   \@namedef{dgo@S}{\lamssource{S}}%
   \def\dggeometry{
      \dg@ZTEMP=\dg@XGRID \multiply\dg@ZTEMP\tw@
      \ifnum\dg@YGRID=\z@ \dg@ZTEMP=\tw@
      \else \divide\dg@ZTEMP\dg@YGRID \fi
      \ifnum\dg@ZTEMP>6\relax \dg@ZTEMP=6\relax \fi
      \ifdg@LATEXGEOM\ifnum\dg@ZTEMP>\f@ur \dg@ZTEMP=\f@ur \fi\fi
      \ifnum\dg@ZTEMP<\@ne \dg@ZTEMP=\@ne \fi
      \unitlength=2sp\relax
      \ifnum\dg@ZTEMP<\tw@
         \advance\dg@ZTEMP\@ne
         \multiply\unitlength\dg@YGRID
      \else
         \multiply\unitlength\dg@XGRID \divide\unitlength\dg@ZTEMP
      \fi
      \dg@XGRID=\dg@ZTEMP \dg@YGRID=\tw@
      \dg@rmcommondiv\tw@\dg@XGRID\dg@YGRID
      \divide\unitlength\dg@YGRID \divide\unitlength\@m
      \dg@LATEXGEOMfalse}
   \def\dgt@nee{\dg@DX=\tw@ \dg@DY=\@ne \dg@SIZE=\tw@}%
   \def\dgt@see{\dg@DX=\tw@ \dg@DY=\m@ne \dg@SIZE=\tw@}%
   \def\dgt@nww{\dg@DX=-\tw@ \dg@DY=\@ne \dg@SIZE=\tw@}%
   \def\dgt@sww{\dg@DX=-\tw@ \dg@DY=\m@ne \dg@SIZE=\tw@}%
   \def\dgt@nnne{\dg@DX=\@ne \dg@DY=\thr@@ \dg@SIZE=\@ne}%
   \def\dgt@nnnw{\dg@DX=\m@ne \dg@DY=\thr@@ \dg@SIZE=\@ne}%
   \def\dgt@sssw{\dg@DX=\m@ne \dg@DY=-\thr@@ \dg@SIZE=\@ne}%
   \def\dgt@ssse{\dg@DX=\@ne \dg@DY=-\thr@@ \dg@SIZE=\@ne}%
   \def\dgt@nnnee{\dg@DX=\tw@ \dg@DY=\thr@@ \dg@SIZE=\tw@}%
   \def\dgt@nnnww{\dg@DX=-\tw@ \dg@DY=\thr@@ \dg@SIZE=\tw@}%
   \def\dgt@sssww{\dg@DX=-\tw@ \dg@DY=-\thr@@ \dg@SIZE=\tw@}%
   \def\dgt@sssee{\dg@DX=\tw@ \dg@DY=-\thr@@ \dg@SIZE=\tw@}%
   \def\dgt@nneee{\dg@DX=\thr@@ \dg@DY=\tw@ \dg@SIZE=\thr@@}%
   \def\dgt@nnwww{\dg@DX=-\thr@@ \dg@DY=\tw@ \dg@SIZE=\thr@@}%
   \def\dgt@sswww{\dg@DX=-\thr@@ \dg@DY=-\tw@ \dg@SIZE=\thr@@}%
   \def\dgt@sseee{\dg@DX=\thr@@ \dg@DY=-\tw@ \dg@SIZE=\thr@@}%
   \def\dgt@neee{\dg@DX=\thr@@ \dg@DY=\@ne \dg@SIZE=\thr@@
      \global\dg@LATEXGEOMtrue}%
   \def\dgt@nwww{\dg@DX=-\thr@@ \dg@DY=\@ne \dg@SIZE=\thr@@
      \global\dg@LATEXGEOMtrue}%
   \def\dgt@swww{\dg@DX=-\thr@@ \dg@DY=\m@ne \dg@SIZE=\thr@@
      \global\dg@LATEXGEOMtrue}%
   \def\dgt@seee{\dg@DX=\thr@@ \dg@DY=\m@ne \dg@SIZE=\thr@@
      \global\dg@LATEXGEOMtrue}%
}
\newcount\dg@HORIZ      \newcount\dg@VERT
\newcount\dg@XLPAD      \newcount\dg@YBPAD
\newcount\dg@XRPAD      \newcount\dg@YTPAD
\newcount\dg@X          \newcount\dg@Y
\newcount\dg@XGRID      \newcount\dg@YGRID
\newcount\dg@SIZE
\newcount\dg@USERSIZE
\newcount\dg@DX         \newcount\dg@DY
\newcount\dg@XLBL       \newcount\dg@YLBL
\newcount\dg@XOFFSET    \newcount\dg@YOFFSET
\newcount\dg@LBLOFF
\newcount\dg@XLBLOFF    \newcount\dg@YLBLOFF
\newcount\dg@LBLPOS
\newcount\dg@XTEMP      \newcount\dg@YTEMP
\newcount\dg@ZTEMP
\newcount\dg@XNODE      \newcount\dg@YNODE
\newcount\dg@XEND       \newcount\dg@YEND
\newcount\dg@WEND       \newcount\dg@HEND
\newcount\dg@COUNT
\newbox\dg@NODEBOX
\newoptenvironment*{diagram}{}[1]{%
   \global\advance\dg@COUNT\@ne \typeout{[diagram \the\dg@COUNT:}%
   \let\node=\dg@node \let\\=\dg@cr \let\arrow=\dg@arrow
   \def\dg@BIGNODE{#1}%
   \ifx\dg@BIGNODE\@empty\else
      \setbox\dg@NODEBOX=\hbox{$\dgeverynode{#1}$}%
      \dg@XTEMP=\wd\dg@NODEBOX
      \dg@YTEMP=\ht\dg@NODEBOX \advance\dg@YTEMP\dp\dg@NODEBOX
      \ifnum\dg@YTEMP=\z@ \dg@ZTEMP=\z@
      \else \dg@ZTEMP=\dg@XTEMP \divide\dg@ZTEMP\dg@YTEMP\fi
      \ifnum\dg@ZTEMP>\dgMAXSQUARE
         \ifnum\dg@ZTEMP<\dgMINDBLSQ \dg@XGRID=\thr@@ \dg@YGRID=\tw@
         \else \dg@XGRID=\tw@ \dg@YGRID=\@ne \fi
      \else \dg@XGRID=\@ne \dg@YGRID=\@ne \fi
      \advance\dg@XTEMP\dgHORIZPAD \advance\dg@YTEMP\dgVERTPAD
      \dg@XLPAD=\dg@XTEMP \divide\dg@XLPAD\tw@
      \dg@XRPAD=\dg@XTEMP \divide\dg@XRPAD\tw@
      \dg@YBPAD=\dg@YTEMP \divide\dg@YBPAD\tw@
      \dg@YTPAD=\dg@YTEMP \divide\dg@YTPAD\tw@
      \advance\dg@XTEMP\dgARROWLENGTH
      \divide\dg@XTEMP\dg@XGRID \divide\dg@XTEMP\@m
      \unitlength=1sp\relax \multiply\unitlength\dg@XTEMP
   \fi
   \typeout{saving...}%
   \dg@X=-\@ne \dg@Y=\z@ \dg@HORIZ=\z@\relax
   \let\dg@SLIST=\@empty
   \let\dg@NLIST=\@empty \let\dg@ALIST=\@empty
   \let\dg@PASS=\dg@savepass
}{
   \dg@VERT=-\dg@Y\relax
   \typeout{calculating...}%
   \ifx\dg@BIGNODE\@empty
      \dg@XGRID=\z@ \dg@YGRID=\z@
      \dg@XLPAD=\z@ \dg@XRPAD=\z@ \dg@YBPAD=\z@ \dg@YTPAD=\z@
      \let\dg@PASS=\dg@geompass
      \dg@NLIST \dg@ALIST
      \ifnum\dg@XGRID>\z@\else \dg@XGRID=\quad\relax \fi
      \ifnum\dg@YGRID>\z@\else \dg@YGRID=\quad\relax \fi
      \dggeometry
      \ifdim\unitlength>\z@\else \unitlength=\quad \fi
      \ifnum\dg@XGRID>\z@\else \dg@XGRID=\@ne \fi
      \ifnum\dg@YGRID>\z@\else \dg@YGRID=\@ne \fi
   \fi
   \dg@LBLOFF=\dgLABELOFFSET \divide\dg@LBLOFF\unitlength
   \multiply\dg@HORIZ\@m \multiply\dg@HORIZ\dg@XGRID
   \multiply\dg@VERT\@m \multiply\dg@VERT\dg@YGRID
   \divide\dg@XLPAD\unitlength \divide\dg@XRPAD\unitlength
   \divide\dg@YBPAD\unitlength \divide\dg@YTPAD\unitlength
   \typeout{drawing...}%
   \let\dg@PASS=\dg@drawpass
   \hspace*{\dg@XLPAD\unitlength}%
   \vcenter{%
      \hsize=0pt\relax\parindent=0pt\relax
      \parskip=0pt\relax\baselineskip=0pt\relax
      \vspace*{\dg@YTPAD\unitlength}%
      \begin{picture}(0,0)(0,0)\dg@NLIST\dg@ALIST\end{picture}%
      \raisebox{\z@}[\z@][\dg@VERT\unitlength]{}%
      \vspace*{\dg@YBPAD\unitlength}%
   }
   \hspace*{\dg@HORIZ\unitlength}%
   \hspace*{\dg@XRPAD\unitlength}%
   \typeout{done]}}%

\def\dg@savepass{s}
\def\dg@geompass{g}
\def\dg@drawpass{d}

\newoptcommand{\dg@node}{\@ne}[2]{%
   \ifx\dg@PASS\dg@savepass
      \dg@XTEMP=#1\relax
      \ifnum\dg@XTEMP<\@ne \dg@XTEMP=\@ne\fi
      \advance\dg@X\dg@XTEMP
      \ifnum\dg@HORIZ<\dg@X \dg@HORIZ=\dg@X \fi
      \setbox\dg@NODEBOX=\hbox{$\dgeverynode{#2}$}%
      \dg@XTEMP=\wd\dg@NODEBOX \advance\dg@XTEMP\dgHORIZPAD
      \dg@YTEMP=\ht\dg@NODEBOX \advance\dg@YTEMP\dp\dg@NODEBOX
      \advance\dg@YTEMP\dgVERTPAD
      \toks\z@=\expandafter{\dg@SLIST}%
      \edef\dg@SLIST{\the\toks\z@
         ,\noexpand\dg@XNODE=\number\dg@X\noexpand\relax
         \noexpand\dg@YNODE=\number\dg@Y\noexpand\relax
         \noexpand\dg@XTEMP=\number\dg@XTEMP\noexpand\relax
         \noexpand\dg@YTEMP=\number\dg@YTEMP\noexpand\relax}%
      \toks\z@=\expandafter{\dg@NLIST}%
      \toks\tw@={\dg@node{#2}}%
      \edef\dg@NLIST{\the\toks\z@
         \noexpand\dg@X=\number\dg@X\noexpand\relax
         \noexpand\dg@Y=\number\dg@Y\noexpand\relax
         \the\toks\tw@}%
   \else\ifx\dg@PASS\dg@geompass
      \ifnum\dg@X=\z@
         \dg@getnodesize
            {\dg@SLIST}{\dg@X}{\dg@Y}{\dg@WEND}{\dg@HEND}%
         \divide\dg@WEND\tw@
         \ifnum\dg@XLPAD<\dg@WEND \dg@XLPAD=\dg@WEND \fi\fi
      \ifnum\dg@X=\dg@HORIZ
         \dg@getnodesize
            {\dg@SLIST}{\dg@X}{\dg@Y}{\dg@WEND}{\dg@HEND}%
         \divide\dg@WEND\tw@
         \ifnum\dg@XRPAD<\dg@WEND \dg@XRPAD=\dg@WEND \fi\fi
      \ifnum\dg@Y=\z@
         \dg@getnodesize
            {\dg@SLIST}{\dg@X}{\dg@Y}{\dg@WEND}{\dg@HEND}%
         \divide\dg@HEND\tw@
         \ifnum\dg@YTPAD<\dg@HEND \dg@YTPAD=\dg@HEND \fi\fi
      \ifnum\dg@Y=-\dg@VERT\relax
         \dg@getnodesize
            {\dg@SLIST}{\dg@X}{\dg@Y}{\dg@WEND}{\dg@HEND}%
         \divide\dg@HEND\tw@
         \ifnum\dg@YBPAD<\dg@HEND \dg@YBPAD=\dg@HEND \fi\fi
   \else\ifx\dg@PASS\dg@drawpass
      \dg@XNODE=\dg@X \multiply\dg@XNODE\@m
      \multiply\dg@XNODE\dg@XGRID
      \dg@YNODE=\dg@Y \multiply\dg@YNODE\@m
      \multiply\dg@YNODE\dg@YGRID
      \setbox\dg@NODEBOX=\hbox{$\dgeverynode{#2}$}%
      \put(\dg@XNODE,\dg@YNODE){\dg@makebox{\box\dg@NODEBOX}}%
   \fi\fi\fi}%

\newoptcommand{\dg@cr}{\@ne}[1]{%
   \ifx\dg@PASS\dg@savepass
      \dg@YTEMP=#1\relax
      \ifnum\dg@YTEMP<\@ne \dg@YTEMP=\@ne \fi
      \advance\dg@Y -\dg@YTEMP\relax
      \dg@X=-\@ne\relax\fi}%
\newoptcommand{\dg@arrow}{\@ne}[2]{%
   \begingroup
   \dg@USERSIZE=#1\relax
   \ifnum\dg@USERSIZE<\@ne \dg@USERSIZE=\@ne \fi
   \dg@parse{#2}%
   \ifx\dg@PASS\dg@savepass
      \ifx\dg@label\dgl@b \let\dg@label=\dgl@t \fi
      \ifx\dg@label\dgl@r \let\dg@label=\dgl@l \fi
      \let\dg@process=\dg@save
   \else\ifx\dg@PASS\dg@geompass
      \let\dg@process=\dg@ignore
      \dg@geomcalc
   \else\ifx\dg@PASS\dg@drawpass
      \let\dg@process=\dg@draw
      \dg@drawcalc
   \fi\fi\fi
   \dg@label{\dg@process{#1}{#2}}}%

\newoptcommand{\arrow}{\@ne}[2]{%
   \dg@parse{#2}%
   \ifx\dg@label\dgl@b \let\dg@label=\dgl@t \fi
   \ifx\dg@label\dgl@r \let\dg@label=\dgl@l \fi
   \dg@label{\dg@textarrow{#1}{#2}}}%

\def\dg@textarrow#1#2#3#4{%
   \mathop{{\dgHORIZPAD=0pt\relax\dgVERTPAD=0pt\relax
      \begin{diagram}
         \node{}\arrow[#1]{#2}{#3}{#4}\node{}
      \end{diagram}}}}

\def\dg@parse#1{%
   \let\dg@label=\dgl@ \dgo@
   \let\dg@type=\@empty \let\dg@lbltype=\@empty
   \@for\dg@list:=#1\do{%
      \ifx\dg@type\@empty \let\dg@type=\dg@list
      \else\ifx\dg@lbltype\@empty \let\dg@lbltype=\dg@list
         \@ifundefined{dgo@\dg@list}{}{\@nameuse{dgo@\dg@list}}%
      \else
         \@ifundefined{dgo@\dg@list}{}{\@nameuse{dgo@\dg@list}}%
      \fi\fi}%
   \@ifundefined{dgt@\dg@type}{\dgt@e}{\@nameuse{dgt@\dg@type}}%
   \@ifundefined{dgl@\dg@lbltype}{}{%
      \dg@letname\dg@label{dgl@\dg@lbltype}}}

\def\dg@draw#1#2#3#4{%
   \put(\dg@X,\dg@Y){\dg@makebox{%
      \begin{picture}(0,0)%
         \thinlines
         \put(\dg@XOFFSET,\dg@YOFFSET){%
            \dg@VECTOR(\dg@DX,\dg@DY){\dg@SIZE}}%
         \put(\dg@XLBL,\dg@YLBL){\dg@makebox{%
            \begin{picture}(0,0)%
               \put(\dg@XLBLOFF,\dg@YLBLOFF){%
                  \dg@makebox[\dg@LBLONE]{$\dgeverylabel{#3}$}}%
               \put(-\dg@XLBLOFF,-\dg@YLBLOFF){%
                  \dg@makebox[\dg@LBLTWO]{$\dgeverylabel{#4}$}}%
            \end{picture}}}%
      \end{picture}}}%
   \endgroup}%

\def\dg@save#1#2#3#4{%
   \endgroup 
   \toks\z@=\expandafter{\dg@ALIST}%
   \toks\tw@={\dg@arrow[#1]{#2}{#3}{#4}}%
   \edef\dg@ALIST{\the\toks\z@%
      \noexpand\dg@X=\number\dg@X\noexpand\relax
      \noexpand\dg@Y=\number\dg@Y\noexpand\relax
      \the\toks\tw@}}%

\def\dg@ignore#1#2#3#4{\endgroup}

\def\dg@geomcalc{%
   \dg@XEND=\dg@SIZE \multiply\dg@XEND\dg@USERSIZE
   \ifnum\dg@DX=\z@
      \dg@YEND=\dg@XEND \dg@XEND=\z@
      \dg@changesign\dg@YEND\dg@DY
   \else
      \dg@changesign\dg@XEND\dg@DX \dg@YEND=\dg@XEND
      \multiply\dg@YEND\dg@DY \divide\dg@YEND\dg@DX
   \fi
   \advance\dg@XEND\dg@X \advance\dg@YEND\dg@Y
   \dg@getnodesize
      {\dg@SLIST}{\dg@XEND}{\dg@YEND}{\dg@WEND}{\dg@HEND}%
   \dg@XOFFSET=\dg@WEND \dg@YOFFSET=\dg@HEND
   \dg@getnodesize
      {\dg@SLIST}{\dg@X}{\dg@Y}{\dg@WEND}{\dg@HEND}%
   \advance\dg@XOFFSET\dg@WEND \divide\dg@XOFFSET\tw@
   \advance\dg@YOFFSET\dg@HEND \divide\dg@YOFFSET\tw@
   \dg@XTEMP=\dgARROWLENGTH \dg@YTEMP=\dgARROWLENGTH
   \ifnum\dg@DX>\z@
      \dg@ZTEMP=\dg@DX \multiply\dg@XTEMP\dg@DX
   \else \dg@ZTEMP=-\dg@DX \multiply\dg@XTEMP -\dg@DX \fi
   \ifnum\dg@DY>\z@
      \advance\dg@ZTEMP\dg@DY \multiply\dg@YTEMP\dg@DY
   \else \advance\dg@ZTEMP -\dg@DY \multiply\dg@YTEMP -\dg@DY\fi
   \ifnum\dg@ZTEMP=\z@\else
      \divide\dg@XTEMP\dg@ZTEMP \divide\dg@YTEMP\dg@ZTEMP
      \advance\dg@XOFFSET\dg@XTEMP \advance\dg@YOFFSET\dg@YTEMP
   \fi
   \divide\dg@XOFFSET\dg@SIZE \divide\dg@XOFFSET\dg@USERSIZE
   \divide\dg@YOFFSET\dg@SIZE \divide\dg@YOFFSET\dg@USERSIZE
   \ifnum\dg@DX=\z@ \dg@XOFFSET=\z@ \fi
   \ifnum\dg@DY=\z@ \dg@YOFFSET=\z@ \fi
   \ifnum\dg@XGRID<\dg@XOFFSET \global\dg@XGRID=\dg@XOFFSET\fi
   \ifnum\dg@YGRID<\dg@YOFFSET \global\dg@YGRID=\dg@YOFFSET\fi
   \relax}

\def\dg@drawcalc{%
   \dg@XEND=\dg@SIZE \multiply\dg@XEND\dg@USERSIZE
   \ifnum\dg@DX=\z@
      \dg@YEND=\dg@XEND \dg@XEND=\z@
      \dg@changesign\dg@YEND\dg@DY
   \else
      \dg@changesign\dg@XEND\dg@DX \dg@YEND=\dg@XEND
      \multiply\dg@YEND\dg@DY \divide\dg@YEND\dg@DX
   \fi
   \advance\dg@XEND\dg@X \advance\dg@YEND\dg@Y
   \dg@getnodesize
      {\dg@SLIST}{\dg@XEND}{\dg@YEND}{\dg@WEND}{\dg@HEND}%
   \divide\dg@WEND\unitlength \divide\dg@HEND\unitlength
   \multiply\dg@DX\dg@XGRID \multiply\dg@DY\dg@YGRID
   \dg@rmcommondiv\tw@\dg@DX\dg@DY
   \dg@rmcommondiv\tw@\dg@DX\dg@DY 
   \dg@rmcommondiv\thr@@\dg@DX\dg@DY
   \multiply\dg@SIZE\dg@USERSIZE \multiply\dg@SIZE\@m
   \ifnum\dg@DX=\z@
      \multiply\dg@SIZE\dg@YGRID
      \divide\dg@HEND\tw@ \advance\dg@SIZE -\dg@HEND
      \dg@getnodesize
         {\dg@SLIST}{\dg@X}{\dg@Y}{\dg@WEND}{\dg@YOFFSET}%
      \divide\dg@YOFFSET\unitlength \divide\dg@YOFFSET\tw@
      \advance\dg@SIZE -\dg@YOFFSET
      \dg@XOFFSET=\z@
      \def\dg@LBLONE{r}\def\dg@LBLTWO{l}%
      \dg@XLBL=\z@ \dg@YLBL=\dg@SIZE
      \multiply\dg@YLBL\dg@LBLPOS
      \divide\dg@YLBL\dgARROWPARTS\relax
      \advance\dg@YLBL\dg@YOFFSET
      \dg@changesign\dg@YLBL\dg@DY
      \dg@changesign\dg@YOFFSET\dg@DY
   \else
      \multiply\dg@SIZE\dg@XGRID
      \ifnum\dg@DY=\z@
         \divide\dg@WEND\tw@ \advance\dg@SIZE -\dg@WEND
         \dg@getnodesize
            {\dg@SLIST}{\dg@X}{\dg@Y}{\dg@XOFFSET}{\dg@HEND}%
         \divide\dg@XOFFSET\unitlength \divide\dg@XOFFSET\tw@
         \advance\dg@SIZE -\dg@XOFFSET
         \dg@YOFFSET=\z@
         \def\dg@LBLONE{b}\def\dg@LBLTWO{t}%
         \dg@YLBL=\z@ \dg@XLBL=\dg@SIZE
         \multiply\dg@XLBL\dg@LBLPOS
         \divide\dg@XLBL\dgARROWPARTS\relax
         \advance\dg@XLBL\dg@XOFFSET
         \dg@changesign\dg@XLBL\dg@DX
         \dg@changesign\dg@XOFFSET\dg@DX
      \else
         \divide\dg@WEND\tw@ \divide\dg@HEND\tw@
         \multiply\dg@HEND\dg@DX \divide\dg@HEND\dg@DY
         \ifnum\dg@HEND<\z@ \multiply\dg@HEND\m@ne \fi
         \ifnum\dg@WEND<\dg@HEND \advance\dg@SIZE -\dg@WEND
         \else \advance\dg@SIZE -\dg@HEND \fi
         \dg@getnodesize
            {\dg@SLIST}{\dg@X}{\dg@Y}{\dg@WEND}{\dg@HEND}%
         \divide\dg@WEND\unitlength \divide\dg@WEND\tw@
         \divide\dg@HEND\unitlength \divide\dg@HEND\tw@
         \multiply\dg@HEND\dg@DX \divide\dg@HEND\dg@DY
         \ifnum\dg@HEND<\z@ \multiply\dg@HEND\m@ne \fi
         \ifnum\dg@WEND<\dg@HEND \dg@XOFFSET=\dg@WEND
         \else \dg@XOFFSET=\dg@HEND \fi
         \advance\dg@SIZE -\dg@XOFFSET
         \dg@changesign\dg@XOFFSET\dg@DX
         \dg@YOFFSET=\dg@XOFFSET
         \multiply\dg@YOFFSET\dg@DY \divide\dg@YOFFSET\dg@DX
         \def\dg@LBLONE{br}\def\dg@LBLTWO{tl}%
         \ifnum\dg@DX<\z@ \ifnum\dg@DY>\z@
            \def\dg@LBLONE{bl}\def\dg@LBLTWO{tr}\fi\fi
         \ifnum\dg@DX>\z@ \ifnum\dg@DY<\z@
            \def\dg@LBLONE{bl}\def\dg@LBLTWO{tr}\fi\fi
         \dg@XLBL=\dg@SIZE
         \multiply\dg@XLBL\dg@LBLPOS
         \divide\dg@XLBL\dgARROWPARTS\relax
         \dg@changesign\dg@XLBL\dg@DX
         \dg@YLBL=\dg@XLBL
         \multiply\dg@YLBL\dg@DY \divide\dg@YLBL\dg@DX
         \advance\dg@XLBL\dg@XOFFSET
         \advance\dg@YLBL\dg@YOFFSET
      \fi
   \fi
   \dg@XLBLOFF=-\dg@DY \dg@changesign\dg@XLBLOFF\dg@DX
   \dg@YLBLOFF=\dg@DX \dg@changesign\dg@YLBLOFF\dg@DX
   \ifnum\dg@DX=\z@ \dg@XLBLOFF=\m@ne \fi
   \dg@XTEMP=\dg@DX \dg@changesign\dg@XTEMP\dg@DX
   \dg@YTEMP=\dg@DY \dg@changesign\dg@YTEMP\dg@DY
   \ifnum\dg@YTEMP>\dg@XTEMP \dg@XTEMP=\dg@YTEMP \fi
   \ifnum\dg@XTEMP=\z@ \dg@XTEMP=\@ne \fi
   \multiply\dg@XLBLOFF\dg@LBLOFF \divide\dg@XLBLOFF\dg@XTEMP
   \multiply\dg@YLBLOFF\dg@LBLOFF \divide\dg@YLBLOFF\dg@XTEMP
   \multiply\dg@X\@m \multiply\dg@X\dg@XGRID
   \multiply\dg@Y\@m \multiply\dg@Y\dg@YGRID
   \relax}%

\def\dg@rmcommondiv#1#2#3{%
   \dg@XTEMP=#2\relax
   \divide\dg@XTEMP #1\relax \multiply\dg@XTEMP #1\relax
   \dg@YTEMP=#3\relax
   \divide\dg@YTEMP #1\relax \multiply\dg@YTEMP #1\relax
   \ifnum\dg@XTEMP=#2\relax \ifnum\dg@YTEMP=#3\relax
      \divide#2#1\relax \divide#3#1\relax \fi\fi}%

\def\dg@changesign#1#2{%
   \ifnum #2<\z@ \multiply#1\m@ne
   \else\ifnum #2=\z@ #1=\z@ \fi\fi}%

\def\dg@getnodesize#1#2#3#4#5{%
   #4=\z@\relax #5=\z@\relax
   \expandafter\@for\expandafter\dg@trynode
   \expandafter:\expandafter=#1\do{%
      \dg@XNODE=\m@ne 
      \dg@trynode
      \ifnum #2=\dg@XNODE \ifnum #3=\dg@YNODE
         #4=\dg@XTEMP\relax #5=\dg@YTEMP\relax\fi\fi}}%

\newoptcommand{\dg@makebox}{}[2]{%
   \expandafter\makebox\expandafter(\expandafter
      0\expandafter,\expandafter0\expandafter)\expandafter
      [#1]{#2}}%

\def\dg@novector(#1,#2)#3{}%

\def\dg@letname#1#2{%
   \relax\expandafter
   \let\expandafter #1\csname #2\endcsname\relax}%

\def\dgl@#1{#1{}{}}%
\def\dgl@t#1#2{#1{#2}{}}%
\def\dgl@b#1#2{#1{}{#2}}%
\def\dgl@tb#1#2#3{#1{#2}{#3}}%
\def\dgl@l#1#2{#1{#2}{}}%
\def\dgl@r#1#2{#1{}{#2}}%
\def\dgl@lr#1#2#3{#1{#2}{#3}}%
\makeatother

\begin{document}
               
%
\expandafter\ifx\csname amssym.def\endcsname\relax \else\endinput\fi
%
\expandafter\edef\csname amssym.def\endcsname{%
       \catcode`\noexpand\@=\the\catcode`\@\space}
\catcode`\@=11
%

\def\undefine#1{\let#1\undefined}
\def\newsymbol#1#2#3#4#5{\let\next@\relax
 \ifnum#2=\@ne\let\next@\msafam@\else
 \ifnum#2=\tw@\let\next@\msbfam@\fi\fi
 \mathchardef#1="#3\next@#4#5}
\def\mathhexbox@#1#2#3{\relax
 \ifmmode\mathpalette{}{\m@th\mathchar"#1#2#3}%
 \else\leavevmode\hbox{$\m@th\mathchar"#1#2#3$}\fi}
\def\hexnumber@#1{\ifcase#1 0\or 1\or 2\or 3\or 4\or 5\or 6\or 7\or
8\or
 9\or A\or B\or C\or D\or E\or F\fi}

\font\tenmsa=msam10
\font\sevenmsa=msam7
\font\fivemsa=msam5
\newfam\msafam
\textfont\msafam=\tenmsa
\scriptfont\msafam=\sevenmsa
\scriptscriptfont\msafam=\fivemsa
\edef\msafam@{\hexnumber@\msafam}
\mathchardef\dabar@"0\msafam@39
\def\dashrightarrow{\mathrel{\dabar@\dabar@\mathchar"0\msafam@4B}}
\def\dashleftarrow{\mathrel{\mathchar"0\msafam@4C\dabar@\dabar@}}
\let\dasharrow\dashrightarrow
\def\ulcorner{\delimiter"4\msafam@70\msafam@70 }
\def\urcorner{\delimiter"5\msafam@71\msafam@71 }
\def\llcorner{\delimiter"4\msafam@78\msafam@78 }
\def\lrcorner{\delimiter"5\msafam@79\msafam@79 }
\def\yen{{\mathhexbox@\msafam@55 }}
\def\checkmark{{\mathhexbox@\msafam@58 }}
\def\circledR{{\mathhexbox@\msafam@72 }}
\def\maltese{{\mathhexbox@\msafam@7A }}

\font\tenmsb=msbm10
\font\sevenmsb=msbm7
\font\fivemsb=msbm5
\newfam\msbfam
\textfont\msbfam=\tenmsb
\scriptfont\msbfam=\sevenmsb
\scriptscriptfont\msbfam=\fivemsb
\edef\msbfam@{\hexnumber@\msbfam}
\def\Bbb#1{{\fam\msbfam\relax#1}}
\def\widehat#1{\setbox\z@\hbox{$\m@th#1$}%
 \ifdim\wd\z@>\tw@ em\mathaccent"0\msbfam@5B{#1}%
 \else\mathaccent"0362{#1}\fi}
\def\widetilde#1{\setbox\z@\hbox{$\m@th#1$}%
 \ifdim\wd\z@>\tw@ em\mathaccent"0\msbfam@5D{#1}%
 \else\mathaccent"0365{#1}\fi}
\font\teneufm=eufm10
\font\seveneufm=eufm7
\font\fiveeufm=eufm5
\newfam\eufmfam
\textfont\eufmfam=\teneufm
\scriptfont\eufmfam=\seveneufm
\scriptscriptfont\eufmfam=\fiveeufm
\def\frak#1{{\fam\eufmfam\relax#1}}
\let\goth\frak

\csname amssym.def\endcsname
\newsymbol\rtimes 226F
\newsymbol\nmid 232D
\newsymbol\varnothing 203F
\newread\epsffilein    
\newif\ifepsffileok    
\newif\ifepsfbbfound   
\newif\ifepsfverbose   
\newdimen\epsfxsize    
\newdimen\epsfysize    
\newdimen\epsftsize    
\newdimen\epsfrsize    
\newdimen\epsftmp      
\newdimen\pspoints     
\pspoints=1bp          
\epsfxsize=0pt         
\epsfysize=0pt         
\def\epsfbox#1{\global\def\epsfllx{72}\global\def\epsflly{72}%
   \global\def\epsfurx{540}\global\def\epsfury{720}%
   \def\lbracket{[}\def\testit{#1}\ifx\testit\lbracket
   \let\next=\epsfgetlitbb\else\let\next=\epsfnormal\fi\next{#1}}%
\def\epsfgetlitbb#1#2 #3 #4 #5]#6{\epsfgrab #2 #3 #4 #5 .\\%
   \epsfsetgraph{#6}}%
\def\epsfnormal#1{\epsfgetbb{#1}\epsfsetgraph{#1}}%
\def\epsfgetbb#1{%
%
%
\openin\epsffilein=#1
\ifeof\epsffilein\errmessage{I couldn't open #1, will ignore it}\else
%
doesn't
%
   {\epsffileoktrue \chardef\other=12
    \def\do##1{\catcode`##1=\other}\dospecials \catcode`\ =10
    \loop
       \read\epsffilein to \epsffileline
       \ifeof\epsffilein\epsffileokfalse\else
%
%
          \expandafter\epsfaux\epsffileline:. \\%
       \fi
   \ifepsffileok\repeat
   \ifepsfbbfound\else
    \ifepsfverbose\message{No bounding box comment in #1; using
defaults}\fi\fi
   }\closein\epsffilein\fi}%
%
%
\def\epsfclipstring{}
\def\epsfclipon{\def\epsfclipstring{ clip}}%
\def\epsfclipoff{\def\epsfclipstring{}}%
\def\epsfsetgraph#1{%
   \epsfrsize=\epsfury\pspoints
   \advance\epsfrsize by-\epsflly\pspoints
   \epsftsize=\epsfurx\pspoints
   \advance\epsftsize by-\epsfllx\pspoints
%
picture.
%
   \epsfxsize\epsfsize\epsftsize\epsfrsize
   \ifnum\epsfxsize=0 \ifnum\epsfysize=0
      \epsfxsize=\epsftsize \epsfysize=\epsfrsize
      \epsfrsize=0pt
%
arithmetic!
reasonably
%
     \else\epsftmp=\epsftsize \divide\epsftmp\epsfrsize
       \epsfxsize=\epsfysize \multiply\epsfxsize\epsftmp
       \multiply\epsftmp\epsfrsize \advance\epsftsize-\epsftmp
       \epsftmp=\epsfysize
       \loop \advance\epsftsize\epsftsize \divide\epsftmp 2
       \ifnum\epsftmp>0
          \ifnum\epsftsize<\epsfrsize\else
             \advance\epsftsize-\epsfrsize \advance\epsfxsize\epsftmp
\fi
       \repeat
       \epsfrsize=0pt
     \fi
   \else \ifnum\epsfysize=0
     \epsftmp=\epsfrsize \divide\epsftmp\epsftsize
     \epsfysize=\epsfxsize \multiply\epsfysize\epsftmp
     \multiply\epsftmp\epsftsize \advance\epsfrsize-\epsftmp
     \epsftmp=\epsfxsize
     \loop \advance\epsfrsize\epsfrsize \divide\epsftmp 2
     \ifnum\epsftmp>0
        \ifnum\epsfrsize<\epsftsize\else
           \advance\epsfrsize-\epsftsize \advance\epsfysize\epsftmp
\fi
     \repeat
     \epsfrsize=0pt
    \else
     \epsfrsize=\epsfysize
    \fi
   \fi
%
parse.
   \ifepsfverbose\message{#1: width=\the\epsfxsize,
height=\the\epsfysize}\fi
   \epsftmp=10\epsfxsize \divide\epsftmp\pspoints
   \vbox to\epsfysize{\vfil\hbox to\epsfxsize{%
      \ifnum\epsfrsize=0\relax
        \includegraphics{#1}%
      \else
        \epsfrsize=10\epsfysize \divide\epsfrsize\pspoints
        \includegraphics{#1}%
      \fi
      \hfil}}%
\global\epsfxsize=0pt\global\epsfysize=0pt}%
%
%
{\catcode`\%=12
\global\let\epsfpercent=
%
%
\long\def\epsfaux#1#2:#3\\{\ifx#1\epsfpercent
   \def\testit{#2}\ifx\testit\epsfbblit
      \epsfgrab #3 . . . \\%
      \epsffileokfalse
      \global\epsfbbfoundtrue
   \fi\else\ifx#1\par\else\epsffileokfalse\fi\fi}%
%
%
\def\epsfempty{}%
\def\epsfgrab #1 #2 #3 #4 #5\\{%
\global\def\epsfllx{#1}\ifx\epsfllx\epsfempty
      \epsfgrab #2 #3 #4 #5 .\\\else
   \global\def\epsflly{#2}%
   \global\def\epsfurx{#3}\global\def\epsfury{#4}\fi}%
%
%
\def\epsfsize#1#2{\epsfxsize}
%
%
\let\epsffile=\epsfbox
  
\addtolength{\baselineskip}{5pt}
\newcommand{\jac}{{\cal J}}
\newcommand{\mod}{{\cal M}_g}
\newcommand{\univ}{{\cal C}_g}
\newcommand{\lra}{\longrightarrow}
\newcommand{\ms}{\mapsto}
\newcommand{\ra}{\rightarrow} 
\newcommand{\pr}{{\bf P}} 
\newcommand{\cl}{{\cal L } }
\newcommand{\siegel }{{\frak h}_g}
\setlength{\unitlength}{1mm}
\newtheorem{lem}{Lemma}[section]
\newtheorem{cor}{Corollary}[section]
\newtheorem{theor}{Theorem}[section]
\newtheorem{prop}{Proposition}[section]
\newtheorem{defi}{Definition}[section]
\newtheorem{rem}{Remark}[section]
\renewcommand{\thetheor}{\Alph{theor}}

\setcounter{section}{-1}

$ $
\vskip.3in
\noindent
\begin{center}
{\cal {Theta line bundles and the determinant of the Hodge bundle.}} \\
\vskip.1in
Alexis Kouvidakis \\
University of Crete \\
Department of Mathematics \\
71409, Iraklion-Crete, Greece\\
e-mail: kouvid@math.uch.gr\\
\end{center}
\vskip.3in

\noindent
Abstract: {\em {In this note we examine the question of expressing the 
determinant of the push forward of a  symmetric line bundle on an abelian 
fibration in terms of the pull back of the relative
dualizing sheaf via the zero section.}}

\section{Introduction.}
Let ${\rm f}: X \lra B$ be a fibration of abelian varieties with a zero 
section $s: B \lra X$. Let $\cl $ be a symmetric line bundle on $X$, 
trivialized along the zero section, which defines a polarization of type 
$D=(d_1, \ldots ,d_g)$  on the fibration. A theorem of Faltings-Chai, see 
Theorem 5.1 in ${\rm Ch.}\,1$ of \cite{FC}, gives the following expression of 
the determinant of the push forward of the line bundle $\cl $ in terms of the 
pull back $s^* \omega _{X/B}$ of the relative dualizing sheaf of the family:
\[ 8 d^3 \det {\rm f}_* \cl = -4 d^4 s^* \omega _{X/B}, 
\]   
where $d:= d_1 \cdots d_g$. 

In this note we show that in the complex analytic category, the above torsion
factor can be improved. More specifically, we have
\begin{theor}
\label{theorem1}
Let ${\rm f} : X \lra B$ be a fibration of complex abelian varieties of 
relative  dimension $g$ and $s$ its zero section. Let  $\cl $ be a symmetric 
line bundle on $X$, trivialized along the zero section, which defines a 
polarization of type $D=(d_1, \ldots ,d_g)$, where $d_1|\ldots |d_g$ are 
positive integers. Let $d=d_1 \cdots d_g$, then 
 \[ 8 \det {\rm f}_* \cl = - 4 d s^* \omega _{X/B}, 
\]  
except when $3|d_g $ and $gcd(3,d_{g-1})=1$, in which case we have that 
  \[ 24  \det {\rm f}_* \cl = -12 d s^* \omega _{X/B}. 
\]
\end{theor}

\noindent   
Moreover, when $\cl $ is totally symmetric i.e. $\cl $ is fiberwise the square 
of a symmetric line bundle, we have 
\begin{theor}
\label{theorem2}
Let ${\rm f} : X \lra B$ be a fibration of complex abelian varieties of 
relative dimension $g \geq 3$ and $s$ its zero section. Let  $\cl $ be a 
totally symmetric line bundle on $X$, trivialized along the zero section, 
which defines a polarization of type $D=(d_1, \ldots ,d_g)$, where 
$d_1|\ldots |d_g$ are positive even integers. Let $d=d_1 \cdots d_g$, then 
\[ \det {\rm f}_* \cl = - \frac{d}{2} s^* \omega _{X/B}, 
\]  
except when $3|d_g $ and $gcd(3,d_{g-1})=1$, in which case we have that 
  \[ 3  \det {\rm f}_* \cl = -3 \frac{d}{2} s^* \omega _{X/B}. 
\]
\end{theor}

\noindent
The theorems are proved by using a ``higher rank" version of the Theta
Transformation Formula, see Propositions \ref{proposition1.1} and 
\ref{proposition1.2}, in order to  construct transition functions for 
$\det {\rm f}_* (\cl )$, see Lemma \ref{lemma2.1}. 

At the end, we apply Theorem \ref{theorem2} to the case of the universal 
Jacobian variety ${\rm f} _{g-1}: \jac ^{g-1} \lra \mod $. This is an abelian 
torsor  which parametrizes line bundles of degree $g-1$ on the fibers of the  
universal curve, where $ \mod $ denotes the moduli space of smooth,  
irreducible  curves of genus $g \geq 3$, without automorphisms.
On that we have a canonical theta divisor defined by the push forward of the
universal symmetric product of degree $g-1$, via the Abel-Jacobi map. We denote
by $\Theta $ the corresponding line bundle and by $\lambda $ the determinant 
of the Hodge bundle of the universal curve $ \psi : {\cal C} \lra \mod $, i.e. 
$ \lambda  =  \det \psi _*\omega _{ {\cal C} / \mod } $.  We then have
\begin{theor}
\label{theorem3}
Keeping the above notation, we have that
\[ \det {\rm f}_{g-1 \; *} (\Theta ^{\otimes n}) =\frac{1}{2}n^g(n-1) \lambda .
\]
\end{theor}     
\noindent
We also give an  alternative way of proving Theorem \ref{theorem3} by 
utilizing special properties of the universal Jacobian varieties. 

\vskip.2in
\noindent
{\bf Acknowledgment.} I would like to thank Professor L. Moret-Bailly for
showing me the alternative way for proving Theorem \ref{theorem3}.

\setcounter{section}{0}

\section{Abelian varieties and theta functions.}
We begin by recalling some standard theory for abelian varieties and theta 
functions. The main reference we follow is the book by Lange and Birkenhake  
\cite{LB}. We also keep (most of) the notation of that book. We denote by 
$X=V/ \Lambda $  an  abelian variety, $V$ is a ${\Bbb C}$-vector space of 
dimension $g$ and  $\Lambda $ a $2g$-full lattice  in $V$.

\subsection{Line bundles on abelian varieties.}    
\label{Line bundles on abelian varieties.}    
On the abelian variety $X=V/\Lambda $, a line bundle is given, up to isomorphy,
by data  $(H,\chi )$, where $H$ is a hermitian  form on $V$ s.t.  ${\rm Im }
H(\Lambda, \Lambda) \subset {\Bbb Z}$ and  $\chi : \Lambda \lra {\Bbb C}_1$  
a semicharacter. The isomorphy class of line bundles on $X$ given by the above 
data is denoted by $L(H , \chi )$, see ${\rm Ch.} \, 2$, $\& \, 2$ in 
\cite{LB}. 
                
Let $\phi : X^{'}=V{'}/ \Lambda{'} \lra X = V / \Lambda $ be a map of abelian 
varieties. To that one can associate two maps, namely the analytic 
representation $\phi _a : V^{'} \lra V$ and the rational representation  
$\phi _r : \Lambda ^{'}  \lra  \Lambda $, see ${\rm Ch.}\, 1$,  $ \& \,2$ in  
\cite{LB}. If $\phi $ is a map  of abelian varieties as above, then  
$\phi ^*L(H, \chi ) =  L(\phi ^*_a H , \phi ^*_r \chi )$.

\subsection{Decomposition of $V$ - Characteristics.}   
\label{Decomposition of $V$.}   
Given a hermitian form $H$ as in Section 
\ref{Line bundles on abelian varieties.}, there is a base $< \lambda _1, 
\ldots , \lambda _g \, ; \, \lambda_{g+1}, \ldots , \lambda_{2g}>$ of 
$\Lambda $, w.r.t. which ${\rm Im}H$ is given by a matrix of the form 
$ \left( \begin{array}{c} \;\; 0 \;\;\;\, D \\ -D \;\;\; 0  \end{array} 
\right)  $, where $D$ is an integral diagonal matrix with elements $d_1| 
\ldots |d_g$. Such a base is called a symplectic base for $H$. When the matrix 
$D$ can be chosen to be the identity, then $H$ is called a principal 
polarization of $X$. 

Let $\Lambda _1$ be the lattice spanned by the first $g$ vectors and 
$\Lambda _2$ the lattice spanned by the last $g$ vectors. If $\lambda \in
\Lambda $, we write $\lambda = \lambda _1 + \lambda _2 $ for the corresponding
decomposition of $\lambda$.  Then $V$, viewed as an ${\Bbb R}$-vector space, 
decomposes as $V=V_1 \oplus V_2 $, where $ V_i = \Lambda _i \otimes {\Bbb R}$; 
we call that one decomposition of $V$ for $H$. If $v\in V$, we write 
$v=v_1 + v_2$ for the corresponding decomposition of $v$. We denote by 
$E:={\rm Im}H$ the alternating form associated to $H$ and by $B$ the symmetric 
form on $V$ which is the ${\Bbb C}$-bilinear extension of the symmetric form 
$H|V_2 \times V_2$, see Lemma 2.1 of ${\rm Ch.} \, 3$ in \cite{LB}. 
We define $\Lambda (H) := \{ v \in V ; {\rm Im} H(v, \Lambda ) \subset 
{\Bbb Z} \}$. Then  $\Lambda(H)= \Lambda(H)_1 \oplus \Lambda(H)_2$, where 
$\Lambda(H)_i :=\Lambda(H) \cap V_i$.

We denote by ${\rm Pic}^H(X)$ the group of line bundles, up to isomorphy, 
corresponding to the Hermitian form $H$. {\em After choosing a decomposition 
of $V$ for $H$ }, we define in  ${\rm Pic}^H(X)$  a distinguished isomorphy 
class of  line bundles, denoted by $L_0 \cong L(H, \chi _0)$, by setting 
$\chi _0 (\lambda ) =  e(\pi i E( \lambda _1, \lambda _2))$. If $L(H, \chi )$ 
is any other element of  ${\rm Pic}^H(X)$, then there exists a $c\in V$, 
unique up to translations by elements of $\Lambda (H)$, s.t. $L(H, \chi ) = 
T^*_c  L(H, \chi _0)$, where $T_c$ is the translation map  by $c$ on $X$. Note 
that if $L(H, \chi )$ is of characteristic $c$ then $\chi (\lambda )= 
\chi _0 (\lambda ) e(2 \pi i E(c, \lambda ))$, see ${\rm Ch.} \, 3$, 
$\& \, 1$ in \cite{LB}.

\subsection{Canonical - Classical factor of automorphy.} 
\label{Canonical - Classical factor of automorphy.} 
A factor of automorphy is a holomorphic map $f :\Lambda \times  V \lra 
{\Bbb C}^{\times }$ satisfying the condition $f(\lambda _1 + \lambda _2 , v) 
= f (\lambda _1, \lambda _2 + v) f( \lambda_2 , v)$. Two factors of automorphy 
$f$ and $f^{'}$ are called  equivalent iff there is a holomorphic function 
$h : V \ra {\Bbb C}^{\times }$ satisfying $f^{'}(\lambda , v)= f(\lambda, v) 
h(v) h(v+\lambda)^{-1}$. We use the notation $f^{'} = f \star h $  to indicate
the above situation. 
           
Let $L(H, \chi )$ be an isomorphy class of line bundles on $X$. We define the 
canonical factor of automorphy $a_{(H, \chi)} : \Lambda \times V  \lra 
{\Bbb C}^{\times }$ by 
\[ a_{(H, \chi)}(\lambda, v):=\chi (\lambda ) \, 
e(\pi H(v, \lambda ) +  \frac{\pi }{2} H(\lambda , \lambda )).  
\] 
We define the classical factor of automorphy $e_{(H, \chi )} : \Lambda \times 
V \lra  {\Bbb C}^{\times }$ by 
\[ e_{(H, \chi )}(\lambda , v):=\chi (\lambda ) \, e (\pi (H-B)(v, \lambda ) + 
\frac {\pi }{2} (H-B)(\lambda , \lambda )). 
\]
We have that 
$ e_{(H, \chi )}(\lambda , v)= a_{(H, \chi)}(\lambda , v) \, e(\frac {\pi }{2} 
B(v,v)) \, e(\frac {\pi }{2} B(v+\lambda , v+ \lambda))^{-1}.$ This shows that 
the above factors are equivalent. When the semicharacter is the $\chi _0$ i.e. 
the characteristic is $0$, we simplify the notation for the above factors of 
automorphy to $a_H$ and $e_H$.

\subsection{Period matrices.} 
\label{Period matrices.} 
Let $X=V/\Lambda $ be an abelian variety and $H$ a polarization of type $D$. 
We choose a symplectic base $< \lambda _1, \ldots , \lambda _g \, ;\, 
\lambda_{g+1}, \ldots , \lambda_{2g}>$ of $\Lambda $ for $H$. Since 
$<\lambda _{g+1}, \ldots , \lambda _{2g}> $ is a ${\Bbb C}$-base of $V$, there 
exists a unique $g \times 2g$ matrix $\Pi $ of the form $\Pi = (Z, D)$, where 
$Z \in \siegel $, with $<\lambda _1 , \ldots , \lambda _g > \, = \, 
<\lambda _{g+1}, \ldots , \lambda _{2g}> Z$. We call $\Pi$ the period matrix 
of $X$ corresponding to the above choice of a symplectic base. Therefore to 
the triple $(X, H , < \lambda _1, \ldots , \lambda _g \, ; \, \lambda_{g+1}, 
\ldots , \lambda_{2g}> )$ corresponds an element $Z \in \siegel $. Conversely, 
given $Z \in \siegel $, let $\Lambda _Z := (Z,D) {\Bbb Z}^{2g}$. Define $X_Z 
:= {\Bbb C}^g / \Lambda _Z$. Let $H_Z$ be the Hermitian form with matrix 
$({\rm Im}Z)^{-1}$ w.r.t. the standard base of ${\Bbb C}^g$.  This one defines 
then a  polarization on $X_Z$ of type $D$.  Let $\lambda _i $ be the vector 
corresponding to the $i$-column of the matrix $(Z,D)$. Then $< \lambda _1, 
\ldots , \lambda _g \, ; \, \lambda_{g+1}, \ldots , \lambda_{2g}>$ turns out 
to be a symplectic base for $H_Z$. Therefore to $Z\in \siegel $ corresponds 
a triple $(X, H_Z , < \lambda _1, \ldots , \lambda _g \, ; \, \lambda_{g+1}, 
\ldots , \lambda_{2g}> )$ and the two constructions are inverse to each other.

\subsection{Canonical - Classical theta functions.} 
\label{Canonical - Classical theta functions.}
Let $f : \Lambda \times V \lra {\Bbb C}^{\times }$ be a factor of automorphy. 
A theta function corresponding to $f$ is a holomorphic function $\theta : 
V \lra {\Bbb C}$ satisfying the functional equation $\theta (\lambda + v ) 
= f(\lambda , v) \theta (v)$. Let $L(H, \chi )$ be of characteristic $c$. 
Theta functions corresponding to the canonical (resp. classical) factor of 
automorphy $a_{(H, \chi)}$ (resp. $e_{(H, \chi )}$) are called canonical 
(resp. classical) theta functions for $L(H, \chi )$. After choosing a
decomposition of $V$ for $H$, we define 
\begin{eqnarray*}
\lefteqn{ \theta ^c_0 (v):= } \\
& & e(-\pi H(v,c) - \frac{\pi }{2} H(c,c) + \frac  {\pi }{2} 
B(v+c, v+c)) \sum _{\lambda _1\in \Lambda _1}e(\pi (H-B)(v+c, \lambda _1) 
- \frac {\pi }{2} (H-B)(\lambda _1, \lambda _1)). 
\end{eqnarray*} 
We have the following facts, see ${\rm Ch.} \, 3$, $\& \, 2$ in \cite{LB}.
\begin{enumerate}
\item  $\theta ^c_0 $ is a canonical theta function for $L(H, \chi)$. By the
definition we deduce that $\theta ^c_0(v)=e(-\pi  H(v,c)-\frac{\pi }{2}H(c,c))
\theta (v+c)$. 
\item Let $\theta ^c_w (v) := a_{(H, \chi)}(w,v)^{-1} 
\,  \theta  ^c_0 (v+w)$, where $\bar{w} \in K(H):=\Lambda(H)/\Lambda $. This 
is well defined and it is also a canonical theta function for $L(H, \chi)$. 
\item The set $ < \theta ^c_w \, ; \, \bar{w} \in K(H)_1:= \Lambda(H)_1 / 
\Lambda _1 > $ forms a base of the canonical theta functions for $L(H, \chi )$. 
\end{enumerate}
 
Let $Z \in \siegel $. We denote by $j_Z$ the isomorphism of ${\Bbb R}$-vector
spaces $j_Z: {\Bbb R}^{2g} \lra {\Bbb C}^g$ which sends $x \ms (Z,1)x$. 
Let $\Lambda _D$ denote the lattice ${\Bbb Z}^g \oplus D {\Bbb Z}^g$ in 
${\Bbb R}^{2g}$ and $\Lambda _Z:=j_Z (\Lambda _D )$ the corresponding lattice 
in ${\Bbb C}^g$. Given $v\in {\Bbb C}^g$, we can write uniquely $v=Zv^1+v^2$, 
with $v^i \in {\Bbb R}^g$. We then have that $\lambda \in \Lambda _Z$ can be 
written  uniquely in the form $\lambda =Z \lambda ^1 + \lambda ^2$ with
$\lambda ^1 \in {\Bbb Z}^g$ and $\lambda ^2 \in D{\Bbb Z}^g$

Let now $L(H, \chi )$, where $H = H_Z$ is a polarization of type $D$, be an
isomorphy class of line bundles of characteristic $c$. We have the following 
facts (most of them can be found in ${\rm Ch.} \, 8$, $\& \, 5$ in \cite{LB}; 
the rest is straight forward calculation): 
\begin{enumerate}
\item $H(v,w) = \, ^tv \, ({\rm Im}Z)^{-1} \,  \bar{w}$, 
      $\;\; B(v,w) = \, ^tv \,  ({\rm Im}Z)^{-1} \, w$. \\ 
      $(H-B)(v,w)=-2i \, ^tv w^1 $, 
      $\;\; E(v,w)=\, ^t v^1 w^2 - \, ^t v^2 w^1 $.
\item $e_{H}(\lambda , v)= e(-\pi i \, ^t \lambda ^1Z\lambda ^1 
      -2\pi i\,^tv \lambda ^1)$. \\
      $a_{H}(\lambda ,v) = e(\pi i \,^t\lambda ^1 \lambda ^2 + 
      \pi \, ^tv({\rm Im}Z)^{-1}\bar{\lambda } + \frac{\pi }{2} 
      \, ^t \lambda ({\rm Im}Z)^{-1} \bar{\lambda } )$. \\
      It is $e_H = a_H \star h$, with 
      $h(v)= e(\frac{\pi }{2} \, ^tv ({\rm Im}Z)^{-1}v)$. 
\item Let ${\Bbb Z}_D$ denote the group ${\Bbb Z}_D := {\Bbb Z}_{d_1} \oplus
      \ldots {\Bbb Z}_{d_g}$. \\
      $\Lambda (H)_1 = \{ Zv^1 \;| \; v^1 \in D^{-1} {\Bbb Z}^g \}$, 
      $\;\; \Lambda (H)_2 = \{ v^2 \;| \; v^2 \in  {\Bbb Z}^g \}.$  \\
      $ K(H)_1 \cong  D^{-1} {\Bbb Z}_D $,
      $\;\; K(H)_2 \cong {\Bbb Z}_D  $. 
\item Let $c = Z c^1 + c^2$. Then $ { \theta} 
      \left[ \begin{array}{c} 
      c^1 \\ c^2 \end{array} \right] (v,Z) := 
      e(-\frac{\pi }{2} B(v,v)+ \pi i \, ^tc^1 c^2 ) 
      \; \theta ^c_0 (v)$ is a classical theta function for $L(H, \chi )$. We
      have 
      \[  \theta \left[ \begin{array}{c} c^1 \\ c^2 \end{array} 
      \right] (v,Z) = \sum _{\lambda ^1 \in {\Bbb Z}^g} e \left(  \pi i 
      \, ^t(\lambda ^1 + c^1)Z(\lambda ^1 +c^1) + 
      2 \pi i \, ^t(v+c^2)(\lambda ^1+c^1) \right) .\]  
\item Let $c=Zc^1+c^2 \in {\Bbb C}^g$ and $w=Zw^1+w^2 \in \Lambda(H)$. 
      Then \\
      a) $ \;\;    \theta  \left[ \begin{array}{c} 
      c^1 + w^1\\ c^2  \end{array} \right] (v,Z) = 
      e(-\frac{\pi }{2} B(v,v)+ \pi i \, ^tc^1 c^2 ) 
      \; \theta ^c_{Zw^1} (v)       $ \\
      b) $ \;\; \theta ^{c+w}(v) = e(- \pi i \, ^tw^1c^2+ \pi i 
         \, ^tc^1w^2 + \pi i \, ^tw^1 w^2 ) \; \theta ^c_w(v)$. \\
      c) $ \;\; \theta ^c_{Zw^1+w^2} (v) =  \theta ^c_{Zw^1}(v)$.  \\
      d) $ \;\; \theta  \left[ \begin{array}{c} 
      c^1 + w^1\\ c^2+w^2  \end{array} \right] (v,Z)=e(2\pi i \, ^tc^1w^2+2
      \pi i \, ^tw^1w^2) \; \theta  \left[ \begin{array}{c} 
      c^1 + w^1\\ c^2  \end{array} \right] (v,Z)$    
\item The set $< \theta \left[ \begin{array}{c} c^1+ D^{-1}m \\ c^2 
      \end{array} \right], \, m \in {\Bbb Z}_D >$ is a base for classical 
      theta functions for $L(H, \chi )$, where $c=Zc^1+c^2$ is the 
      characteristic of $L(H, \chi )$. 
\end{enumerate}

\subsection{Action of the symplectic group.}  
\label{Action of the symplectic group.}  
Let $D$ denote a polarization type. Let $\Gamma _D := \{ R \in M_{2g}({\Bbb Z})
, \; R  \left( \begin{array}{c} \;\; 0 \;\;\;\; D \\ -D  \;\;\; 0 \end{array}  
\right) \, ^tR = \left( \begin{array}{c} \;\; 0 \;\;\;\;  D \\ -D  \;\;\; 0 
\end{array}  \right) \}$.
Let $G_D :=  \{ M \in {\rm Sp}_{2g}({\Bbb Q}), \; M= \left( \begin{array}{c} 
I \;\;\; 0 \\ 0  \;\;\; D \end{array}  \right) ^{-1} R \left( \begin{array}{c}
I \;\;\; 0 \\ 0  \;\;\; D  \end{array}  \right), \; \mbox{for some} \;\; R 
\in \Gamma_D \}$. Note that for  $R= \left( \begin{array}{c} A \;\;\; B \\ 
\Gamma  \;\;\; \Delta \end{array}     \right) \in \Gamma _D$ and $M = \left( 
\begin{array}{c} \alpha \;\;\; \beta \\ \gamma  \;\;\; \delta  \end{array}
\right) \in G_D$ defined as above , we have 
\[ \alpha =A , \;\; \beta =BD, \;\; \gamma = D^{-1} \Gamma , \;\; \delta =
D^{-1} \Delta D.
\]

The group $G_D$ is characterized as  $G_D =\{ M \in {\rm Sp}_{2g}({\Bbb Q}), 
\; ^tM \Lambda _D \subset \Lambda _D  \} $. The group $G_D$  acts on 
$\siegel $ as follows, see ${\rm Ch.} \, 8$, $\& \, 1$ in \cite{LB}. Let
$M = \left( \begin{array}{c} \alpha \;\;\; \beta \\ \gamma  \;\;\; \delta 
\end{array}     \right) \in G_D$; we define $M(Z):=(\alpha Z + \beta D)
(\gamma Z + \delta )^{-1}$. Two abelian varieties $X_Z$ and $X_{Z^{'}}$ of 
type $D$ are isomorphic iff the matrices $Z$ and $Z^{'}$ lie in the same orbit 
under the above action i.e. $Z^{'}=M(Z)$. The isomorphism is given by $\phi 
(M) : X_{Z} \lra X_{M(Z)}$ s.t. the corresponding map $\phi (M)_r : 
\Lambda _{Z} \lra \Lambda _{M(Z)}$ has matrix $R= \, ^tM^{-1} $ w.r.t. the 
canonical symplectic bases defined by $Z$ and $M(Z)$.  It turns out that the 
map $\phi (M)_a : {\Bbb C}^g \lra {\Bbb C}^g$ has matrix $A^{-1}$, where 
$A= \, ^t(\gamma Z + \delta )$, w.r.t. the standard base of ${\Bbb C}^g$. In 
addition, $\phi (M) _a^*H_{M(Z)} = H_{Z}$. We therefore have the following 
diagram
\noindent
\begin{equation}
\begin{diagram}[{\cal C}_{g,d}]
\node{{\Bbb R}^{2g} } \arrow{e,t}{j_Z } \arrow{s,l}{\phi(M)_r } 
\node{{\Bbb C}^g} \arrow{s,r}{\phi (M)_a } \arrow{e}
\node{X_Z} \arrow{s,r}{\phi (M)}  \\
\node{{\Bbb R}^{2g}} \arrow{e,t}{j_{M(Z)}} 
\node{{\Bbb C}^g} \arrow{e}
\node{X_{M(Z)}}
\end{diagram}
\end{equation}
\noindent     

Define $M_Z(v):= A^{-1}v\, (=\phi (M)_a(v))$. We then have that the action of 
$G_D$ on $\siegel $, lifts to an action on ${\Bbb C}^g \times \siegel $ by 
defining $M (v, Z):=(M_Z(v), M(Z))$. Over the space $\siegel $ there exists 
a universal family ${\frak X}_D \lra \siegel $ with fiber over $Z$ to be the 
abelian variety $X_Z$. This is defined as follows. $\Lambda _D = {\Bbb Z}^g 
\oplus D {\Bbb Z}^g$ acts on ${\Bbb C}^g \times \siegel $ by $n(v,Z)=(v + 
j_Z(n) ,Z)$. This action is free and properly discontinuous and the quotient 
space is ${\frak X}_D$.  The above action of $G_D$  on ${\Bbb C}^g \times 
\siegel $, descents to an action on ${\frak X}_D$.

\subsection{Factors of automorphy and line bundles - Sections and theta
functions.}
\label{Factors of automorphy and line bundles.}
Let $\pi : V \lra X$ denote  the canonical map. A factor of automorphy 
$f : \Lambda \times V \ra {\Bbb C}^{\times }$ gives rise to a line bundle 
on $X$ as follows: Take the action of $\Lambda $ on $V \times {\Bbb C}$ 
defined by $\lambda \; (v,z) := (v+\lambda, f(\lambda,v)z)$. We then define 
$L$ to be the quotient space of $V \times {\Bbb C}$ under this action. 

We can also find transition functions for this line bundle:  Take a covering
$X=\cup _{a} U^a$ of $X$ such that the preimage of each $U^a$ under the map
$\pi $ to be a disjoint union of sets isomorphic to $U^a$. For each $a$, we
choose one of those lifts, say $W^a$, and let $\pi _a : W^a \stackrel{\cong }
{\ra }  U^a$ be the restriction of the map $\pi $. Define $w_a(x) := 
\pi _a ^{-1}(x) \in W_a$. For any $a,b $ with $U^{ab} := U^a \cap U^b \neq 
\emptyset $, there is a unique $\lambda _{ab} \in \Lambda $ satisfying 
$\lambda _{ab} := w_b(x) -  w_a(x)$ for all $x \in U^{ab}$. In that way we 
get  an identification  of $W^a \cap \pi ^{-1} (U^{ab})$ with  $W^b \cap 
\pi ^{-1}(U^{ab})$. Define now holomorphic functions $g_{ab} : U^{ab} \ra  
{\Bbb C}^{\times }$ given by $g_{ab}(x) :=  f(\lambda _{ab}, w_a(x))$, which 
satisfy the cocycle condition for transition functions. Let $L^{tr} :=  
\bigsqcup _a U^a \times {\Bbb C} / \sim $, where $(x_a,z_a) \sim (x_b,z_b)$ 
iff $x_a=x=x_b$ and $z_b=g_{ab}(x)z_a$, be the standard construction of a 
line bundle by its transition functions. There is an isomorphism of $L^{tr}$  
with $L$ given by $L^{tr} \ni (x_a ,z) \ms [w_a(x_a) , z] \in L$ with inverse 
$L \ni [v ,z] \ms (x_a=\pi (v) , f(v-w_a(x_a),w_a(x_a))^{-1} z) \in L^{tr}$.
In that way the $g_{ab}$'s become transition functions of $L$ w.r.t. the 
trivialization induced by composing the natural trivializations of $L^{tr}$
with the above defined inverse map.  Note that by making different choices of 
the liftings $W^a$ of $U^a$ in the above, we get equivalent transition 
functions.      

Let $f$ and $f^{'}$ be equivalent factors of automorphy  satisfying 
$f^{'}(\lambda , v)= f(\lambda, v) h(v) h(v+\lambda)^{-1}$, where $h : V \lra 
{\Bbb C}^{\times }$ holomorphic. By fixing the above covering data we define 
$h^a : U^a \ra {\Bbb C}^{\times }$ by $h_a(x_a) :=h(w_a(x_a))$. The line 
bundles $L$ and $L^{'}$ corresponding to the above factors of automorphy are 
then isomorphic and their transition functions (constructed as above w.r.t. 
the common fixed covering data) satisfy $g^{'}_{ab}(x)=g_{ab}(x)
\frac{h_a(x)}{h_b(x)}$. There is also a canonical choice of isomorphism 
$\Phi _h : L \lra L^{'}$: By viewing the line bundles as quotient spaces, 
via the action of the lattice, the isomorphism is defined by $[v,z] \ms 
[v, h(v)^{-1}z]$. By viewing the line bundles as quotient spaces, via their 
transition functions and the above identification, the isomorphism is defined 
by sending $(x_a, z_a) \ms (x_a, h_a(x_a)^{-1} z_a)$.  

Let $\phi : X^{'} \lra X$ be a map of abelian varieties. If $f$ is a factor of
automorphy for $X$, we denote by $\phi ^* f$ the factor of automorphy for  
$X^{'}$ given by $(\phi _r \times \phi _a)^* f$. If $L$ is the line bundle 
on $X$ corresponding to $f$ then $\phi ^* L $ is the line bunlde on $X^{'}$
corresponding to $\phi ^* f$.  If $\theta $ is a theta function for the factor
of automorphy $f$ then $\phi _a^* \theta $ (or $\phi ^* \theta $ in a
simplified notation) is a theta function for the factor of automorphy 
$\phi ^*f$. 

Sections of the line bundle $L$ given by the factor of automorphy $f$,
correspond to theta functions $\theta : V \ra {\Bbb C}$ satisfying  the
functional equation $\theta (\lambda +v) = f(\lambda,v) \theta (v)$. 
The relation is the following. Given a section $s$ of $L$, let $s(x)=[v_s(x),
z_s(x)]$ (we view $L$ as a quotient space given by the action of the lattice
corresponding to $f$). We then define $\theta _s (v):= f(v-v_s(x), v_s(x))\,
z_s(x) $. It turns out that this is a well defined theta function for the
factor of automorphy $f$. Conversely, given a theta function $\theta $ for
$f$, we define $s(x):=[v(x), \theta (v(x))]$, where $v(x)$ an arbitrary vector
which lies over $x$. This is also well defined and the two processes are
inverse to each other.

If $f^{'}=f \star h$ as in Section 
\ref{Factors of automorphy and line bundles.} and $\Phi _h : L \lra L^{'}$ 
the canonical choice of isomorphism as in Section
\ref{Factors of automorphy and line bundles.} above, then
given a section $s^{'} \in H^0(X,L^{'})$ corresponding to  the theta function 
$\theta _{s^{'}}$, we have that $s:=\Phi _h^{*} s^{'} \in H^0(X,L)$
corresponds to the theta function $\theta _s := h(v) \theta _{f^{'}} (v)$.

\subsection{Theta transformation formulas.}
\label{Theta transformation formulas.}
Let $X_Z$ be an abelian variety with a  line bundle $L(H_Z,\chi )$ which 
defines a polarization $H_Z$ on $X_Z$ . We denote by $c$ the characteristic 
of $L$. Let $M \in G_D$ as in Section 
\ref{Action of the symplectic group.} and define $Z^{'}=M(Z)$. Let 
 $\psi =\psi (M)  : X_{Z^{'}} \lra X_Z $ be the inverse of the map $\phi=
\phi (M) : X_Z \lra X_{Z^{'}}$ defined in  Section 
\ref{Action of the symplectic group.}. 
Then the line bundle $\psi ^* L(H_Z, \chi)$ has type $\psi ^*H_Z = H_{Z^{'}}$,
semicharacter $\chi ^{'} = \psi ^* \chi $ and characteristic $M[c]$ with
$M[c]^1=\delta c^1 - \gamma c^2 + \frac{1}{2}(D\gamma \, ^t\delta)_0$ and
$M[c]^2=-\beta c^1 + \alpha c^2 + \frac{1}{2}(\alpha \, ^t\beta )_0$, see
Lemma 4.1 of Ch.$\,$II in \cite{LB}, where there is an unfortunate omission 
of $D$ in the expression of $M[c]^1$ (The $(\;\;)_0$ stands for the diagonal
vector of the matrix).   

\begin{lem}
\label{lemma1.1} 
We have that $\psi ^* e_{(H_Z, \chi )}=e_{(H_{Z^{'}}, \chi ^{'})} \star
h_{\psi }$, where $ h_{\psi }(v)=e(\pi i \, ^tv ( \gamma Z + \delta )^{-1} 
\gamma  v)$. Also, $\phi ^*e_{(H_{Z^{'}}, \chi ^{'})} = e_{(H_Z, \chi )}  
\star h$, where $h(v)=e(-\pi i \, ^tv (\gamma Z + \delta )^{-1} \gamma  v)$.
\end{lem}

\noindent
{\em Proof.} We have that $a_{(H_Z, \chi )} = e_{(H_Z, \chi )} \star h_1 $,
where $h_1(v)=e(- \frac{\pi }{2} \, ^tv ({\rm Im}Z)^{-1}v)$, see Section 
\ref{Canonical - Classical theta functions.}.
By applying $\psi ^*$ we get that $\psi ^* a _{(H_{Z}, \chi )} =\psi ^*
e_{(H_Z, \chi )}  \star \psi ^*h_1$. Since $\psi ^* a_ {(H_Z, \chi )} = 
a_{(H_{Z^{'}}, \chi ^{'})}$  and $a_{(H_{Z^{'}}, \chi ^{'})} = 
e _{(H_{Z^{'}}, \chi ^{'})} \star h_1 ^{'}$, where $h_1^{'}(v^{'})=
e(-\frac {\pi }{2}\, ^tv^{'} ({\rm Im}Z)^{-1}v^{'})$, we get that  
$h_{\psi }(v)=\psi ^* h_1(v)^{-1} \, h_1^{'}(v^{'})$ i.e. 
$h_{\psi }(v)=e(\frac{ \pi}{2} \, ^tv ({\rm Im}Z)^{-1}v - \frac{\pi }{2} 
\, ^tv^{'} ({\rm Im}Z^{'})^{-1}v^{'})$, where $v^{'}:=\phi _a(v)$. 
A formula given in the proof of Proposition 6.1, ${\rm Ch.} \, 8$ in 
\cite{LB} (or a straight forward calculation),  implies the above form for
$h_{\psi }$. To prove the second formula, we apply $\phi ^*$ to the first one. 
\begin{flushright}
${\Box}$
\end{flushright}
    
The tuple $B^Z:=< \,  \theta ^c_{ZD^{-1}m}(v) \, ; \, m \in {\Bbb Z}_D >$ 
forms a base of the canonical theta functions for $L(H_Z, \chi  )$ and the 
tuple $B^{Z^{'}}:=< \,  \theta ^{M[c]}_{Z^{'}D^{-1}n}(v^{'}) \, ; \, n \in 
{\Bbb Z}_D >$ forms a base of the canonical theta functions for $L(H_{Z^{'}}, 
\chi ^{'} )$. Since $\psi ^*L(H_Z, \chi ) = L(H_{Z^{'}}, \chi ^{'})$ and 
$\psi ^* a _{(H_Z, \chi )} = a_{(H_{Z^{'}}, \chi ^{'})} $,  we get
that $\psi ^*B^Z:=<\psi _a^*  \theta ^c_{D^{-1}m}(v)\,; \, m \in {\Bbb Z}_D > $ 
forms also a base of the canonical theta functions for $L(H_{Z^{'}}, 
\chi ^{'})$. There is then an invertible complex matrix $C$ relating  
$\psi ^*B^Z$  with the base  $B^{Z^{'}}$.

\begin{prop}
\label{proposition1.1}
Keeping the above notation, let we assume that the characteristic 
$c \in \frac{1}{2} \Lambda (H_Z)$. We then have that the matrix $C$ with 
$\psi ^* B^Z = C\, B^{Z^{'}}$   is of the form $C=(\det (\gamma Z + \delta ))
^{-\frac{1}{2}} \; C_M $, where $C_M$ is a matrix depending on $M$ with 
$\det C_M = \zeta _8$, except when $3|d_g$ and $(d_{g-1},3)=1$, in which case 
we have  $\det C_M = \zeta _{24}$ (here by  $\zeta _m $ we denote an $m$-root 
of unity).     
\end{prop}

\noindent
{\em Proof.} Let $G_D^{{\rm int}} := G_D \cap {\rm Sp}_{2g}({\Bbb Z})$. A
matrix $M$ belongs to $G_D^{{\rm int}}$ if $M =  \left( \begin{array}{c} 
I \;\;\; 0 \\ 0  \;\;\; D \end{array}  \right) ^{-1} R \left( \begin{array}{c} 
I \;\;\; 0 \\ 0  \;\;\; D  \end{array}  \right)$, where 
$R =   \left( \begin{array}{c}  A \;\;\; B \\ \Gamma   \;\;\; \Delta 
\end{array}  \right) \in \Gamma _D  $ with $\Gamma = D \Gamma _1 $, $\Gamma _1
\in M_g({\Bbb Z})$.  Let therefore define $\Gamma _D^{\rm int} := \{ R =   
\left( \begin{array}{c}  A \;\;\; B \\ \Gamma   \;\;\; \Delta \end{array}  
\right) \in \Gamma _D  \;\; \mbox{with} \;\; \Gamma = D \Gamma _1 ,\; 
\Gamma _1 \in M_g({\Bbb Z}) \} $. We have the following lemma:

\begin{lem}
\label{lemma1.2}
The group $\Gamma_D$ is generated by $\Gamma _D^{{\rm int}}$ and the matrix 
$J:= \left( \begin{array}{c}  \;0 \;\;\; -I \\ I  \;\; \;\; \;\;0 \end{array}  
\right)$ and therefore  the group  $G_D$ is generated by $G_D^{{\rm int}}$ 
and the matrix $\left( \begin{array}{c} \;\; 0 \;\;\;\; -D \\ D^{-1} \;\;\;\; 
0 \end{array}  \right)$.
\end{lem}

\noindent
{\em Proof.} For the proof we use results from \cite{B}. We have seen that 
$K(D)= D^{-1} {\Bbb Z}_D \oplus {\Bbb Z}_D$. A matrix $R \in \Gamma _D$ 
defines an action on $K(D)$ given by multiplication by $\left( \begin{array}{c} 
I \;\;\; 0 \\ 0  \;\;\; D \end{array}  \right) \, ^t R \left( \begin{array}{c} 
I \;\;\; 0 \\ 0  \;\;\; D  \end{array} \right)^{-1}$. By using the 
identification $K(D) \cong  {\Bbb Z}_D \oplus {\Bbb Z}_D$ via the isomorphism 
$ D^{-1} {\Bbb Z}_D \oplus  {\Bbb Z}_D   \lra {\Bbb Z}_D \oplus  {\Bbb Z}_D $ 
given by the matrix $\left( \begin{array}{c} D \;\;\; 0 \\ 0  \;\;\; I 
\end{array}  \right)$,  the action of a matrix $R \in \Gamma _D$ on $K(D) 
\cong {\Bbb Z}_D \oplus  {\Bbb Z}_D$ is given by multiplication by $\bar{D} 
\, ^t R \bar{D}^{-1}$, where $\bar{D}:= \left( \begin{array}{c} D \;\;\; 0 \\ 
0  \;\;\; D \end{array}  \right)$. One can define on $K(D)$ an alternating 
form $e^D$, see \cite{LB}, \cite{B}, and the above action becomes a symplectic 
action with respect to that form. We denote the symplectic group (w.r.t. $e^D$)
of $K(D) \cong {\Bbb Z}_D \oplus  {\Bbb Z}_D$ by ${\rm Sp}(D)$. We then have 
an exact sequence 
\[ 0 \lra \Gamma _D(D) \lra \Gamma _D \stackrel{\pi }{\lra }{\rm Sp}(D) \lra 0, 
\]  
where $\pi (R) := \bar{D} \, ^t R \bar{D}^{-1}$ and $\Gamma _D(D):= \{ R \in
\Gamma _D \, | \; R=I + \bar{D}A, \, A \in M_{2g}({\Bbb Z}) \}$ (that the
sequence is exact on the right it is proven in \cite{B}). Note that $\Gamma _D
(D) \subset \Gamma _D^{\rm int}$. Therefore it is enough to show that every 
element of ${\rm Sp}(D)$ has a lift to an element of $\Gamma _D$ which is a 
product of the matrix $J$ and elements of  $\Gamma _D^{{\rm int}}$.  

Following the notation of \cite{B}, we have that $A \in L_D$ iff $\bar{D} \, 
^t A \bar{D}^{-1} \in \Gamma_D $, where $L_D$ is defined in Section 2 
of \cite{B}. A matrix $A \in L_D$ acts on $K(D)$ by multiplication. In the 
proof of Theorem 2 of \cite{B}, it is shown that a matrix  $A \in {\rm Sp}(D)$ 
has a lift $\tilde{A} \in L_D $   which satisfies a relation of the form 
\[ \left( \begin{array}{c}  I \;\;\; 0 \\ c_1  \;\; I \end{array}   \right)
\left( \begin{array}{c} \;\; a \;\;\;\;\;\; Dy \\ -D  \;\;\; d \end{array} 
\right)
\left( \begin{array}{c}  U_1 \;\;\; 0 \\ \; 0  \;\;\;\; U_2 \end{array} \right)
\left( \begin{array}{c}  I \;\;\; b_1 \\ 0  \;\;\; I \end{array} \right)
\tilde{A}
= \left( \begin{array}{c}  I \;\;\; b_2 \\ 0  \;\;\; I \end{array} \right),
\]
where $y$ is an integral diagonal matrix and all the matrices belong 
to $L_D$. The inverse  of a matrix $\left( \begin{array}{c}  a \;\;\; b \\ c 
\;\;\; d \end{array}   \right) $ in $L_D$ belongs to $L_D$ and is given by 
$\bar{D} \left( \begin{array}{c}  \;\;\; ^td \;\;\; - \,^t b \\ - \, ^t c  
\;\;\;\;\; ^ta  \end{array}   \right) \bar{D}^{-1} $. Therefore $A$ has a lift 
$R$  in $\Gamma _D$ which is given by $R=\bar{D} \, ^t\tilde{A} \bar{D}^{-1}$. 
By observing now that $\left( \begin{array}{c}  I \;\;\; 0 \\ c  \;\;\; I 
\end{array} \right) = J \left( \begin{array}{c}  -I \;\;\;\;\; c \\ \;\; 0  
\;\;\; -I \end{array}  \right) J$, we conclude that $R=  (JN_1J) N_2 N_3 N_4 
(JN_5J)$, where $N_i \in \Gamma _D^{\rm int}$. This proves the lemma. 
\begin{flushright}
${\Box}$
\end{flushright}  

As in \cite{Mu}, we can rewrite the formula we want to  prove in the form 
\begin{eqnarray*}
\lefteqn{ < \theta ^c_{Z D^{-1}m }(v, Z) \, ; \, m \in  {\Bbb Z}_D > \, 
\sqrt{dv_1 \wedge \ldots \wedge dv_g} = } \\
& & = C_M <\theta ^{M[c]}_{Z^{'} D^{-1}n} (v^{'}, Z^{'}) \, ; \, n \in 
{\Bbb Z}_D > \, \sqrt {dv_1^{'} \wedge \ldots \wedge dv_g^{'}},
\end{eqnarray*}
where $v^{'}= \, ^t(\gamma Z + \delta )^{-1}v$ and $Z^{'}=(\alpha Z+ \beta)
(\gamma Z +\delta )^{-1}$. Note that if $c \in \frac{1}{2} \Lambda(H_Z)$ then 
$M[c] \in \frac{1}{2} \Lambda(H_{Z^{'}})$. We conclude that if the formula 
holds for $M_1, \, M_2 \in G_D$ then it also holds for $M_1M_2$. It is enough 
therefore to verify the formula for matrices in the generators. We express the 
formula we want to prove in terms of classical theta functions by using 
Lemma \ref{lemma1.1} and the formulas in Section 
\ref{Canonical - Classical theta functions.}. It becomes
\begin{eqnarray}
\label{eq:ttf-1}
\lefteqn{ < \theta  \left[  \begin{array}{c} c^1+ D^{-1}m \\ c^2 \end{array}
\right] (v, Z) \, ; \, m \in  {\Bbb Z}_D > = e( -\pi i  \, ^tv (\gamma Z + 
\delta )^{-1} \gamma  v) } \\  
& & e(-\pi i \, ^tM[c]^1M[c]^2 + \pi i \, ^tc^1c^2) (\det (\gamma Z +\delta ))
^{-\frac{1}{2} }  \; C_M  <  \theta  \left[ \begin{array}{c} M[c]^1+D^{-1}n \\ 
M[c]^2  \end{array}  \right] (v^{'},Z^{'}) \, ; \, n \in  {\Bbb Z}_D  >. 
\nonumber 
\end{eqnarray}
  
\vskip.1in
\noindent
{\bf Matrices of the form} $  \left( \begin{array}{l} 0 \;\;\;\;\;\; -D 
\\ D^{-1}  \;\;\;\;\; 0  \end{array}  \right)   $.  
In this case $e( - \pi i  \, ^tv (\gamma Z + \delta )^{-1} \gamma v) = e( -\pi 
i \, ^tv Z^{-1}  v )$ and $ \det (\gamma Z + \delta ) = \frac {\det  Z}{d} $.
We also have that $v^{'}=D Z^{-1}v$, $Z^{'}=- DZ^{-1} D$ and $M[c]^1 = 
-D^{-1}c^2$, $M[c]^2=Dc^1$. The relation (\ref{eq:ttf-1}) in this case becomes 
\begin{eqnarray*}
\lefteqn{ < \theta  \left[  \begin{array}{c} c^1+ D^{-1}m \\ c^2 \end{array}
\right]  (v, Z)  \, ; \, m \in  {\Bbb Z}_D > = e( -\pi i \, ^tv Z^{-1} v) } \\ 
& & e( 2\pi i \, ^tc^1c^2) ( \frac{\det Z}{d} )^{-\frac{1}{2} }  \;  C_M        
< \theta  \left[ \begin{array}{c} -D^{-1}c^2 + D^{-1}n  \\ Dc^1 \end{array} 
\right] (v^{'},Z^{'}) \, ; \, n \in  {\Bbb Z}_D  >.          \nonumber 
\end{eqnarray*}

As in \cite{Mu}, we apply Fourier transform. Write  $ \, \theta  \left[  
\begin{array}{c}  c^1+D^{-1}m  \\ c^2 \end{array}  \right]  (v, Z)= 
\sum _{\lambda \in {\Bbb Z}^g} f(\lambda )  $  where $f(x):= e \left( \pi i \, 
^t(x + c^1+D^{-1}m) Z (x +c^1+ D^{-1}m) + 2 \pi i \, ^t(v+c^2) 
(x +c^1+ D^{-1}m) \right)$. Let $\hat{f}(x):= \int _{{\Bbb R}^g} f(x) \, 
e(2 \pi i \, ^tx \lambda ) \, dx$. We then have that $  \theta  \left[  
\begin{array}{c}  c^1+D^{-1}m \\ c^2 \end{array}  \right]  (v, Z)= 
\sum _{\lambda \in {\Bbb Z}^g } \hat{f}^(\lambda )$. Using Lemma 5.8 of 
${\rm Ch.} \, {\rm II}$ in \cite{Mu}, we get by substituting $x^{'} =  
(x + c^1+D^{-1}m)$ that \[ \hat{f}( \lambda ) = e(- 2 \pi i \, ^t(c^1+D^{-1}m)
\lambda ) \; (\det \frac{Z}{i})^{- \frac{1}{2}} \; e (- \pi i \, ^t(v+c^2 + 
\lambda) Z^{-1} (v+c^2+\lambda)).
\]
We therefore get that
\begin{eqnarray*}
\lefteqn{  \theta  \left[ \begin{array}{c} c^1+D^{-1}m \\ c^2 \end{array}  
\right]  (v, Z)= e(-\pi i  \, ^t(v+c^2) Z^{-1} (v+c^2) ) \, \, 
(\det \frac{Z}{i})^{-\frac{1}{2}}   }   \\ 
& & \sum _{\lambda \in {\Bbb Z}^g} e(-2 \pi i \, ^t (c^1+ D^{-1}m) \lambda  
- 2 \pi i \, ^t(v+c^2) Z^{-1} \lambda - \pi i \; ^t \lambda Z^{-1} \lambda ).
\end{eqnarray*}  
Observe now that by substituting $ - \lambda = D k + n$ with $n \in 
{\Bbb Z}_D$, we can rewrite the sum as 
\begin{eqnarray*}      
\lefteqn{ \sum _{\lambda \in {\Bbb Z}^g} e (- 2 \pi i  \, ^t (c^1+ D^{-1}m)  
\lambda  - 2 \pi i \, ^t(v+c^2) Z^{-1} \lambda - \pi i \; ^t \lambda Z^{-1} 
\lambda ) = }  \\ 
& &   \sum _{n \in  {\Bbb Z}_D}  \sum _{k \in {\Bbb Z}^g} e \left( 2 \pi i \,
^t (c^1+ D^{-1}m) (Dk + n )  + 2 \pi i \, ^t(v+c^2) Z^{-1} (Dk +n)- \pi i \, 
^t (Dk + n)Z^{-1} (Dk + n) \right).
\end{eqnarray*}
Using straight forward calculation, one finally deduces that
\begin{eqnarray*}
\lefteqn{  \theta  \left[ \begin{array}{c} c^1+D^{-1}m \\ c^2 \end{array}  
\right] (v, Z)= e( - \pi i  \,  ^t v Z^{-1} v) \, (\det \frac{Z}{i})
^{-\frac{1}{2}} \, e(2\pi i \, ^tc^1c^2)  } \\
& & \sum _{n \in {\Bbb Z}_D} e(2 \pi i \, ^tm D^{-1}n)\,  \theta  \left[  
\begin{array}{c} -D^{-1}c^2 + D^{-1}n \\ Dc^1 \end{array}  \right]  
(v^{'}, Z^{'}).                                                
\end{eqnarray*}  

The matrix $C_M$ we are asking for has $m,n$ entry equal to $(\frac{d}{i^g})
^{-\frac{1}{2}} \,  e(2\pi i \, ^tmD^{-1}n )$. Let $d:= \det D$. The matrix 
$A=(a_{mn})_{m, n \in {\Bbb Z}_D}$ with $a_{mn}:= e(2 \pi i \, ^tmD^{-1}n)$ 
has determinant $\det A= \zeta _4 d^{\frac {d}{2}}$. One way to see this is 
the  following. Denote by ${\Bbb C}[{\Bbb Z}_{d_i}]$ the ${\Bbb C}$-vector 
space of  dimension $d_i$ ``corresponding" to the finite group 
${\Bbb Z}_{d_i}$. Fix the natural base $<m,\, m \in {\Bbb Z}_{d_i}>$ on that 
vector space. Define the map  $\phi _i : {\Bbb C}[{\Bbb Z}_{d_i}] \lra 
{\Bbb C}[{\Bbb Z}_{d_i}]$  by  $\phi _i (m ):= \sum _{n \in {\Bbb Z}_{d_i}} 
e(2 \pi i n d_i ^{-1} m)n $. Let $C _i$ be the matrix corresponding to 
$\phi _i $. Then $\det C _i= \zeta _4 d_i ^{\frac{d_i}{2}}$, where $\zeta _4$ 
is á $4$-th root of unity. Indeed, it is easy to see that $\det (C_i^2)=
\pm d_i^{d_i}$. Observe now that the matrix $A$ is the matrix corresponding 
to the tensor product of the maps $\phi _i$ and so, its determinant is 
$\det A = \det C_1^{\frac{d}{d_1}} \cdots  \det C_g^{\frac{d}{d_g}} = \zeta _4
d^{\frac {d}{2}}$. 

To conclude this case, observe that $ (\frac{d}{i^g})^{-\frac{d}{2}} =
\zeta _8 d ^{-\frac{d}{2}}$. Therefore $\det C_M =\zeta _8$.  
                            
\vskip.2in
\noindent
{\bf Matrices in $G_D^{{\rm int}}$. }  Let $M\in G^{\rm int}_D$. Then $M$ 
corresponds to an isomorphism $\psi : X_{Z^{'}}/  \Lambda _{Z^{'}} \lra X_Z 
/\Lambda _Z $ which is a lift of an isomorphism of ppav. In this case the usual
theta transformation formula holds, see ${\rm Ch.}\, 8$,  $ \& \,6$ in 
\cite{LB}. Let $ a = c+Zw^1 $, $c \in \frac{1}{2} \Lambda(H), \; w^1=D^{-1}w_1
\in D^{-1}{\Bbb Z}_D $. 

We have the following, see ${\rm Ch.}\, 8$, $ \& \, 4$ and $ \& \,6$ 
in \cite{LB}.
\begin{enumerate}
\item $ \psi ^* \theta ^{a}(v,Z) = C(Z,M,a) \; \theta ^{M[a]}(v^{'},Z^{'})$.
\item $C(Z,M,a) = C(Z,M,0) \; e(\pi i E(M[0], A^{-1}a)) $. 
\item $ C(Z,M,0)=k(M) \; e(\pi i \, ^tM[0]^1M[0]^2) \, 
      \det (\gamma Z + \delta)^{-\frac{1}{2}}$, where $k(M)$ is an 8-th root 
      of unity.  
\end{enumerate}
Note that $M[a]^1=M[c]^1+\delta w^1$ and $M[a]^2=M[c]^2-\beta w^1$. Using the 
above formulas and the formulas in Section 
\ref{Canonical - Classical theta functions.} 
we get 
\begin{eqnarray*}
\lefteqn{  e(-\pi i \, ^tc^2w^1) \psi ^* \theta ^c _{Zw^1}(v,Z) = C(Z,M,0) 
\, e(\pi i E(M[0], A^{-1}a)) } \\ 
& & e(-\pi \, ^t (\delta w^1)M[c]^2+ \pi i \, ^tM[c]^1(-\beta w^1) +
\pi i \, ^t (\delta w^1)(-\beta w^1)) \; \theta^{M[c]}_{Z^{'}\delta w^1} 
(v^{'}, Z^{'}). 
\end{eqnarray*}
We have
\begin{eqnarray*}
 M[0]^1 = \frac{1}{2}(D\gamma \, ^t \delta )_0, \;\;\;\;\;\;\;\;\;\;\;
\;\;\;\;\;\;  
&  & \;\;  M[0]^2 = \frac{1}{2}( \alpha \, ^t \beta )_0 \;\;\;      \\
 (A^{-1}a)^1 = \delta (c^1+w^1)-\gamma c^2, \;\;\;  
&  &  \;\; (A^{-1}a)^2 = - \beta (c^1 + w^1)+ac^2  \\
M[c]^1 =\delta c^1 - \gamma c^2+ \frac{D}{2}( \gamma \, ^t \delta )_0,  
&  &  \;\; M[c]^2 = - \beta c^1 + ac^2+ \frac{1}{2}( \alpha \, ^t \beta )_0.  
\end{eqnarray*}
We then get 
\[ \psi ^* \theta ^c _{ZD^{-1}w_1}(v,Z) = k(M) e(\pi i k) (\pi i \lambda w^1) 
e(-\pi i \, ^t w^1 \beta \delta w^1) \det (\gamma Z + \delta)^{-\frac{1}{2}} 
\theta^{M[c]} _{ Z^{'}D^{-1}\Delta w_1 },
\]
where 
$ k= \, ^t M[0]^1 M[0]^2 + \, ^tM[0]^1 (-\beta c^1+ac^2) - \, ^t M[0]^2 
(\delta c^1 - \gamma c^2) $ and $\lambda =-\,^tM[0]^1\beta -\, ^tM[0]^2 
\delta - \, ^tM[c]^2 \delta -\,  ^tM[c]^1 \beta + \, ^tc^2 $. 
Observe now that $k \in \frac{1}{4d_g} {\Bbb Z}^g$ and $\lambda \in
\frac{1}{2}{\Bbb Z}^g $. 
 
Note that when $\gamma \in M_g ({\Bbb Z})$ i.e. $\Gamma = D \Gamma _1$
for some integral matrix $\Gamma _1$, we have that $\Delta $ acts as a 
permutation on ${\Bbb Z}_D$. Indeed, the relation $\Delta D \,^tA - \Gamma D 
\, ^tB =D$ implies $\Delta (D \, ^tA D^{-1} ) = I + \Gamma  (D \, ^tB D^{-1})$ 
i.e. $\Delta (D \, ^tA D^{-1} ) = I + D \Gamma _1  (D \, ^tB D^{-1} )$ and so, 
$  \Delta A_1 = I + DM $ for some integral matrices $A_1, M$. This implies that
$\Delta $ induces an epimorphism on the group ${\Bbb Z}_D$ and therefore an
automorphism. 

Therefore the matrix $C_M$ of the proposition is going to have in the $w_1,
\Delta w_1$-entry the value $ k(M) e(\pi i k) e(\pi i \lambda D^{-1}w_1)
e(-\pi i \, ^t w_1 \, ^t \Delta D^{-1} B w_1 )$ and any other entry zero.  
We have 
\[  \prod _{w_1 \in {\Bbb Z}_D} e(\pi i \lambda D^{-1}w_1) = e(\pi i 
\sum _{i=1}^g \frac{\lambda _i }{d_i} \sum  _{w_1 \in {\Bbb Z}_D} w_1^i ) =
e(\pi i \sum _{i=1}^g \frac{\lambda _i }{d_i} \frac{d}{d_i} 
\frac{d_i(d_i-1)}{2}) 
\]
The sum which appears above belongs to $\frac{1}{4} {\Bbb Z}$ and so, the
product is a $\zeta _8$. Also, the matrix $\, ^t \Delta D^{-1}B $ is 
symmetric and let $\alpha _{ij} = \frac{a_{ij}}{d_i}, a_{ij} \in {\Bbb Z}$ 
its $ij$-entry. Then  
\[ \prod _{w_1 \in {\Bbb Z}_D} e(\pi i \, ^t w_1 \, ^t \Delta D^{-1} B c_1) =
e \left( \pi i  \sum _{i=1}^{g} \frac{a_{ii}}{d_i} d \frac{(d_i-1)(2d_i-1)}{6} 
+ 2 \sum _{1 \leq i < j \leq g} \frac{a_{ij}}{d_i} d \frac{(d_i-1)(d_j-1)}{4}
\right)
\] 
The sums which appear above belong to $\frac{1}{2} {\Bbb Z}$ and so, the
product is a $\zeta _4$, except when $3|d_g$ and  $(d_{g-1},3)=1$, in which 
case it belongs to $\frac{1}{6} {\Bbb Z}$ and the product is a $\zeta _{12}$.
To conclude, we have  that $\det C _M = \zeta _8 $, except when
$ 3|d_g$ and  $(d_{g-1},3)=1$, in which case $\det C= \zeta _{24}$ 
\begin{flushright}
${\Box}$
\end{flushright}
     
For the case now of a totally symmetric bunlde, note first that in Lemma
\ref{lemma1.1}, if $L(H_Z, \chi _0 ) $  is of characteristic $c=0$,  then 
$\psi ^*L(H_Z , \chi _0) $ is again of characteristic $0$. Indeed, in the 
case of an ``even" polarization, we have always that $\chi _0 =1$ and so,  
$\psi ^*_r \chi _0 = 1 = \chi ^{'} _0$. 

\begin{prop}
\label{proposition1.2}
Keeping the notation of Proposition \ref{proposition1.1}, we assume in
addition that  the characteristic $c=0$, $\cl $ is totally symmetric and 
that $g \geq 3$. We then have that the matrix $C$ with  $\psi ^* B^Z = C\, 
B^{Z^{'}}$   is of the form $C=(\det (\gamma Z + \delta  ))^{-\frac{1}{2}} \; 
C_M $, where $C_M$ is a matrix depending on $M$ with $\det C_M = 1$, except when
$3|d_g$ and $(d_{g-1},3)=1$ in which case we have  $\det C_M = \zeta _{3}$.
\end{prop}

\noindent
{\em Proof.} The proof is a modification of the proof of Proposition
\ref{proposition1.1}. 

At the end of the paragraph ``{\bf Matrices of the form} 
$  \left( \begin{array}{l} 0 \;\;\;\;\;\; -D \\ D^{-1}  \;\;\;\;\; 0  
\end{array}     \right)   $": For $g \geq 3$ the number $\frac{d}{d_i}$ is a 
multiple of $4$ and so, $\det C _i^{\frac{d}{d_i}}=d_i^{\frac{d_i}{2}}$.  
Therefore $\det A = d^{\frac{d}{2}}$.  Also, for $g \geq 3$ we have 
$ (\frac{d}{i^g})^{-\frac{d}{2}} = d ^{-\frac{d}{2}}$  Therefore, 
$\det C_M =1$.   

At the end of the paragraph ``{\bf Matrices in $G_D^{{\rm int}}$ }": The first
sum is an even integer and so, the first product is $1$. For the second sum,  
if $ g \geq 3$ the right summand is an even integer. For the left summand we
have that it is an even integer except when $3|d_g$ and  $(d_{g-1},3)=1$ in 
which case it belongs to $\frac{2}{3} {\Bbb Z}$. Therefore for the second
product we have that it is $1$, except when $3|d_g$ and  $(d_{g-1},3)=1$, in 
which case it is $\zeta _3$.  We also claim that  ${\rm sgn}(\Delta )=1$, 
where by ${\rm sgn}(\Delta )$ we denote the sign of the permutation of 
${\Bbb Z}_D$ defined by the action of $\Delta $.

Indeed, let $d_i = 2^{k_i} m_i$ with $1 \leq k_1 \leq k_2 \leq \ldots 
\leq k_g$ and $m_1 | m_2 | \ldots |m_g$ odd integers. Define  
${\Bbb Z}_{\rm ev}: = {\Bbb Z}_{2^{k_1}} \oplus \cdots \oplus 
{\Bbb Z}_{2^{k_g}} $ a group of order $n_{\rm ev} = 2^{k_1 + \cdots k_g}$ and 
${\Bbb Z}_{\rm odd}:= {\Bbb Z}_{m_1} \oplus \cdots \oplus {\Bbb Z}_{m_g}$ a 
group of order $n_{\rm odd} = m_1 \cdots m_g$. Then ${\Bbb Z}_D = 
{\Bbb Z}_{\rm ev} \oplus {\Bbb Z}_{\rm odd}$. Let $\phi : {\Bbb Z}_D \lra 
{\Bbb Z}_D$ an automorphism. Then $\phi ({\Bbb Z}_{\rm ev}) = 
{\Bbb Z}_{\rm ev}$ and $\phi ({\Bbb Z}_{\rm odd}) = {\Bbb Z}_{\rm odd}$. 
Indeed, ${\Bbb Z}_{\rm ev} = \{ x \in {\Bbb Z}_D \;\; \mbox{with} \;\; 
2^{k_g}x=0 \}$ and ${\Bbb Z}_{\rm odd} = \{ x \in {\Bbb Z}_D \;\; \mbox{with} 
\;\; m_g x=0 \}$. Both those relations are invariant under $\phi $. Let 
therefore $\phi _{\rm ev} $ be the restriction of $\phi $ to 
${\Bbb Z}_{\rm ev}$ and  $\phi _{\rm odd} $ the restriction of $\phi $ to 
${\Bbb Z}_{\rm odd}$. If we interpret $\phi $ as a linear automorphism of the 
vector space  ${\Bbb C}[{\Bbb Z}_D]$ then $\phi = \phi _{\rm ev} \otimes 
\phi _{\rm odd}$. By using the formula for the determinant of the tensor 
product of operators, we get ${\rm sgn}\phi = {\rm sgn}\phi _{\rm ev}
^{n_{\rm odd}} \; {\rm sgn} \phi _{\rm odd}^{n_{\rm ev}}$. Therefore, since 
$n_{ev}$ is an even number, it is enough to prove the result for 
${\Bbb Z}_D = {\Bbb Z}_{\rm ev}$. 

Let now $E$ be the matrix which corresponds to the automorphism $\phi _{\rm
ev}$. We call  elementary transformations of ${\Bbb Z}_{\rm ev}$
those which correspond to left or right multiplication by a matrix
of one of the following forms: $1$ in the diagonal and $a_{ij} \in {\Bbb Z}$ 
in some $ij$-entry if $j\geq i$ (and zero every where else) or $1$ in the 
diagonal and $2^{k_i-k_j} a_{ij}$, $a_{ij} \in {\Bbb Z}$, in some $ij$-entry 
if $j > i$ (and zero every where else). Note that both those types of matrices 
induce automorphisms of ${\Bbb Z}_{\rm ev}$. We then claim that by multiplying 
the matrix $E=(e_{ij})$ with the above type of matrices we can derive a matrix 
with all the elements of the last row, except the diagonal one, to be zero 
mod$2^{k_g}$ and the $i$-th element of the last column, with $1 \leq i<g$, 
to be zero mod$2^{k_i}$. Indeed, we can first assume that $e_{gg}$ is an odd 
integer:  the determinant of $E$ must be an odd number in order $E$ to define 
an automorphism and so, some of the elements of the last row must be odd. If 
$e_{gg}$ is even, let $e_{gj_0}, \, j_0 <g $ be the  odd element. But then 
using an elementary transformation  we can  add the $j_0$-th column to the 
last column and so the $gg$-entry becomes odd. Since $D^{-1}ED$ is an integral 
matrix we have that $e_{gj}=2^{k_g - k_j} m_{gj}, \, m_{gj} \in {\Bbb Z}$. But 
now the equation $2^{k_g - k_j} e_{gg}x \equiv -2^{k_g - k_j} m_{gj} \, 
{\rm mod} 2^{k_g}$ has a solution and therefore by multipling the matrix $E$ 
on the right by the elementary matrix with $2^{k_g - k_j}x$ in the $gj$-entry, 
we get that the $gj$-entry of the product is zero mod$2^{k_g}$. By a similar
argument, by multiplying on the left by an elementary matrix with $x$ in the
$ig$-entry, where $x$ is the solution of $ c_{gg}x \equiv  -c_{ig} \, 
{\rm mod} 2^{k_i}$, we get that the $ig$-entry of the product is zero 
mod$2^{k_i}$. 

A matrix like the one we produced, corresponds to an even permutation of 
${\Bbb Z}_{\rm ev}$. Indeed, by writing ${\Bbb Z}_{\rm ev}$ as a direct 
sum of two groups, the second of which is the ${\Bbb Z}_{2^{k_g}}$, the 
action is a direct sum of actions. Then, by applying the formula for the 
signatures stated above, we get that the signature of the permutation is one, 
since both groups are of even order. A similar argument gives that the action 
which is given by the elementary matrices induces an even permutation (here 
we have to use the hypothesis that $g \geq 3$). We therefore get that the 
permutation given by $\phi _{\rm ev}$ is an even one. This concludes the 
proof of the proposition.
\begin{flushright}
${\Box}$
\end{flushright}  

\setcounter{section}{1}

\section{Abelian fibrations.}
Everything we have stated which holds for a fixed abelian variety $X=V /
\Lambda $, can be transferred easily to hold for a fibration $X \lra U$ of 
abelian varieties of type $D$, with the base $U$ to be a simply connected  
Stein manifold (such as $\siegel $). In this case, the universal covering 
$\tilde{X}$ of $X$  will take the place of $V $ and the homotopy group 
$\pi _1(X)$ will take the place of $\Lambda $. For the universal abelian 
fibration ${\rm f_{un}}:{\frak X}_D \lra \siegel$ we have by construction, see
Section \ref{Action of the symplectic group.}, that $\pi _1 ({\frak X}_D) 
= \Lambda _D $ and $\tilde{{\frak X}_D} = {\Bbb C}^{2g} \times \siegel $.     
On  ${\frak X}_D$  we have a line bundle corresponding to the (holomorphic)
classical factor of automorphy $e: \Lambda _D \times ({\Bbb C}^g \times 
\siegel) \lra {\Bbb C}^{\times }$ given by $e(l;v,Z)=e(-\pi \, i \, ^t 
\lambda ^1 Z \lambda ^1 - 2\pi \, i \, ^tv \lambda ^1)$, where $\lambda ^1$
denotes the first $g$ components of $l$. We denote by $\cl _{\frak X}$ the 
above line bundle. 

\subsection{Canonical form of line bundles on abelian fibrations.}
\label{Canonical form of line bundles on abelian fibrations.}    
Let ${\rm f} : X \lra B$ be a fibration  of abelian varieties  with a 
symmetric line bundle $\cl $ on $X$, trivialized along the zero section, 
which defines a polarization of type $D$ on each fiber. We denote by $s: B 
\lra X$ its zero section and let $S:=s(B)$. We cover $S$ by open analytic 
(connected) sets $V^a$ of $X$. The sets $U^a:=s^*V^a$ define an open covering 
of $B$. Let $(U^a, \lambda ^a_1, \ldots , \lambda ^a_{2g})$ be a choice of 
a symplectic base on the fibers of the restriction of  the abelian fibration 
over $U^a$. For each $a,b$ with $U^{ab} \neq \emptyset $, there is a matrix 
$M^{ab}=\left( \begin{array}{c} \alpha ^{ab} \;\; \beta ^{ab} \\ \gamma ^{ab}
\;\; \delta ^{ab} \end{array} \right) \in G_D$ relating the two symplectic 
bases as in Section \ref{Action of the symplectic group.}. Let $X^a$ be the 
preimage of $U^a$ under ${\rm f}$. We then have a period map $p_a : U^a \lra 
\siegel $ such that the diagram
\noindent
\begin{equation}
\begin{diagram}[{\Bbb C}^g  ]
\node{\tilde {X}^a } \arrow{e,t}{\tilde {\pi  } _a } \arrow{s} 
\node{{\Bbb C}^g \times \siegel} \arrow{s}  \\    
\node{X^a } \arrow{e,t}{\pi   _a } \arrow{s,l}{{\rm f} } 
\node{{\frak X}_D} \arrow{s,r}{{\rm f_{un}} } \\    
\node{U^a} \arrow{e,t}{p _a } \node{\siegel } 
\end{diagram}
\end{equation}
\noindent     
is a diagram of fiber products, where $\tilde{X}^a$ denotes the universal 
cover of $X^a$. We put $Z^a(s):=p_a(s), \; s \in U^a$. Having chosen the 
period map $p_a$, we can identify $\tilde{X}^a$ with ${\Bbb C}^g \times  
p_a(U^a)$ and $\pi _1(X^a)$ with $\Lambda _D$.  We write $(\pi _1(X^a), \, 
\tilde{X}^a) \stackrel{p_a}{\cong} (\Lambda _D , \,{\Bbb C}^g \times  
p_a(U^a))$ to denote  this identification.  We now ``complete" $V^a$ to an 
open cover $X^a=\cup _i V^a_i$  of $X^a$ s.t. $V^a_i \cap S = \emptyset $  
for all $V^a_i \neq V^a$ and such that each $V^a_i$ has an isomorphic lift 
$W^a_i$ in the universal cover $\tilde {X}^a \stackrel{p_a}{\cong} {\Bbb C}^g 
\times  p_a(U^a)$. We can also assume that the intersection of any two of 
those open sets is simply connected.   

Consider the induced diagram
\noindent
\begin{equation}
\begin{diagram}[{\cal C}_{g,d}]
\node{X^{ab} } \arrow{e,t}{\pi _a } \arrow{s,l}{{\rm f} } 
\node{\pi _a(X^{ab})} \arrow{s,l}{{\rm f_{un}} } \arrow{e,t}{\phi _{ab}}
\node{\pi _b(X^{ab})} \arrow{s,r}{{\rm f_{un} }}
\node{X^{ab}} \arrow{w,t}{\pi _b} \arrow{s,r}{{\rm f}}  \\
\node{U^{ab}} \arrow{e,t}{p _a} 
\node{p_a(U^{ab})} \arrow{e,t}{\mu _{ab}}
\node{p_b(U^{ab})}
\node{U^{ab}} \arrow{w,t}{p_b}
\end{diagram}
\end{equation}
\noindent     
We denote by $\mu _{ab}$ the map which sends $Z \ms M^{ab}(Z)$ and by 
$\phi _{ab} $ the map induced by the action of $M^{ab}$ on the universal abelian
variety ${\frak X}$ (fiberwise is the map $\phi (M)$ of 
Section \ref{Action of the symplectic group.}). 
We have that $\pi _b = \phi _{ab} \, \pi _a $ and $p_b = \mu _{ab} \, p_a$. 

Note that a symmetric line bundle $L(H, \chi )$ has always characteristic 
$c \in \frac{1}{2} \Lambda (H)$ w.r.t. any decomposition of $H$. Indeed, that
$L(H, \chi )$ is symmetric is equivalent to $\chi (\lambda )= \pm 1$ for all
$\lambda \in \Lambda $. The claim is a consequence of the definition of
$\chi (\lambda )$. The restriction of  $\cl $ on the fibers over $U^a$ is then
of characteristic $c_a =Z^a(s)c^1_a+c^2_a$, where $(c^1_a,c^2_a) \in 
\frac{1}{2}( D^{-1}{\Bbb Z}^g \oplus {\Bbb Z}^g) $,  w.r.t. the decomposition 
defined by the choice of  the symplectic base. Note that  when $\cl $ is 
totally symmetric the characteristic is going to be $0$.   

Define the line bundles $\cl _a :=\pi ^*_a T^*_{c_a} \cl _{\frak X}$ on  
$X^a$. Note that in the totally symmetric case, we have that $\cl _a :=
\pi ^*_a \cl _{\frak X}$. Then $\cl _a \cong \cl |_{X^a}$ because the 
restriction on the fibers have the same type and characteristic and $U^a$ is 
an open analytic set. Via the identification $(\pi _1(X^a), \, \tilde{X}^a) 
\stackrel{p_a}{\cong} (\Lambda _D , \,{\Bbb C}^g \times  p_a(U^a))$, 
the line bundle  $\cl _a $  corresponds to the factor of automorphy 
$e_a :  \Lambda _D \times ({\Bbb C}^g \times U^a) \lra {\Bbb C}^{\times } $,
where $e_a=e_a(\lambda ^a, v^a ,Z^a(s))$ is the classical factor of 
automorphy of characteristic $c_a$.    

On $X^{ab}$ we have $\cl |_{X^{ab}} \cong \pi _b^*T^*_{c_b} \cl _{\frak X}
\cong \pi ^*_a \phi ^*_{ab} (T^*_{c_b} \cl _{\frak X}|_{\pi _b(X^{ab})})$. 
On the other hand, $\cl |_{X^{ab}} \cong \pi _a^*T^*_{c_a} \cl _{\frak X}$. We
therefore have that $ \phi ^*_{ab} (T^*_{c_b} \cl _{\frak X}|_{\pi _b(X^{ab})})
\cong  T^*_{c_a} \cl _{\frak X}|_{\pi _a(X^{ab})}$. Note that when we fix the 
identification  $(\pi _1(X^{ab}), \, \tilde{X}^{ab}) \stackrel{p_a}{\cong} 
(\Lambda _D , \,{\Bbb C}^g \times  p_a(U^{ab}))$, the line bundle $\cl _a 
|_{X^{ab}}$ corresponds to the factor of automorphy $e_a$ while the line 
bundle $\cl _b|_{X^{ab}}$ corresponds to the factor of automorphy 
$\phi _{ab}^* e_b$. 

By  Lemma \ref{lemma1.1}, we have that  $ \phi ^*_{ab} e_b |_{\Lambda _D 
\times {\Bbb C}^g \times p_a(U^{ab})} = e_a |_{\Lambda _D \times {\Bbb C}^g 
\times p_a(U^{ab})} \star h^{ab}$. In other words  we have that
\[ \phi^*_{ab}e _b=e _a \star h^{ab} \;\; \mbox{on} \;\; (\pi _1(X^{ab}), \, 
\tilde{X}^{ab}) \stackrel{p_a}{\cong} (\Lambda _D , \,{\Bbb C}^g \times  
p_a(U^{ab})),
\] 
where $ h^{ab}(v^a,Z^a(s))= e(- \pi i \, ^tv^a(\gamma ^{ab}Z^a(s) + \delta
^{ab})^{-1} \gamma ^{ab} v^a ). $  Note that $v^b= \, ^t (\gamma Z^a(s) +
\delta)^{-1}v^a$ and $Z^b(s)=(\alpha Z^a(s)+\beta)(\gamma Z^a(s)+\delta)^{-1}$.

We glue now all the line bundles $(\cl _a, X^a)$ together 
to get on $X$ a line bundle $\cl _{\rm can}$ which restricts to 
the trivial bundle on the zero section of the fibration, in the following way: 

Fix a trivialization $g^a_{ij}=e_a(\lambda ^a_{ij},w^a_i,Z^a(s))$  of 
$\cl _a$ on $X^a$ (resp. $g^b_{kl}=e_b(\lambda ^b_{kl},w^b_k,Z^b(s))$  
of $\cl _b$ on $X^b$), by choosing lifts $W^a_i$ of the open sets $V^a_i$ in 
the universal cover $\tilde {X}^a \stackrel{p_a}{\cong} {\Bbb C}^g  \times  
p_a(U^a)$ (resp. lifts $W^b_k$ of $V^b_k$ in the universal cover 
$\tilde {X}^b \stackrel{p_b}{\cong} {\Bbb C}^g \times  p_b(U^b)\,$). 
We then have a canonical isomorphism $\cl _a \stackrel {\simeq }{\lra} 
\bigsqcup _i (V^a_i  \times {\Bbb C})/ \stackrel{a}{\sim } \, $  
(resp. $\cl _b \stackrel {\simeq }{\lra} \bigsqcup _k (V^b_k  \times 
{\Bbb C}) / \stackrel{b}{\sim } \,$), see Section 
\ref{Factors of automorphy and line bundles.}.

Let now $\Phi _{ab} : \cl _a|_{X^{ab}} \lra  \cl _a|_{X^{ab}} $ be the 
canonical isomorphism defined in Section 
\ref{Factors of automorphy and line bundles.} corresponding to 
$ \phi^*_{ab}e_b=e_a \star h^{ab}$. We construct an  explicit isomorphism 
$\psi ^{ab}$ which makes the following diagram commutative:   
\noindent
\begin{equation}
\label{diag:coclboaf-3}
\begin{diagram}[aaaaaaaaaaaaaaa]
\node{\cl _a|_{X^{ab}} } \arrow{e,t}{\Phi _{ab}} \arrow{s,r}{\simeq} 
\node{\cl _b|_{X^{ab}}  } \arrow{s,r}{\simeq}  \\    
\node{\bigsqcup _i (V^a_i  \times {\Bbb C})/ \stackrel{a}{\sim } |_{X^{ab}} } 
\arrow{e,t}{\psi ^{ab}  } 
\node{\bigsqcup _k (V^b_k \times {\Bbb C})/  \stackrel{b}{\sim }|_{X^{ab}}  }   
\end{diagram}
\end{equation}
\noindent  
Such an isomorphism is given by $(x^a_i(s)=x,\, z^a_i) \ms (x^b_k(s)=x, 
\, z^b_k= \psi ^{ab}_{ik}(x) z^a_i)$, where the $\psi ^{ab}_{ik}$'s satisfy the
relation  $\psi ^{ab}_{jl}(x) \, e_a (\lambda ^a_{ij}, w^a_i(x),Z^a(s))
=e_b(\lambda ^b_{kl},w^b_k(x),Z^b(s)) \, \psi^{ab}_{ik}(x)$ and the cocycle 
conditions $\psi ^{ab}_{ik} \, \psi^{ba}_{ki}=1, \; \psi^{ab}_{ik} \, 
\psi^{bc}_{ks} \, \psi^{ca}_{si}=1$ (we follow here the notation of Section 
\ref{Factors of automorphy and line bundles.}).

Let $\pi ^a_i : W^a_i \stackrel{\simeq}{\lra } V^a_i$ denote the canonical
maps. For each $V^{ab}_{ik}:=V^a_i \cap V^b_k \neq \emptyset$,  we have that 
$\pi_{i}^{a\, -1}(V^{ab}_{ik})\subset W^a_i$ is a lift of $V^{ab}_{ik}$ in 
the universal cover $\tilde{X}^{ab}\stackrel{p_a}{\cong} {\Bbb C}^g \times 
p_a(U^{ab})$. We also have that $\tilde{\phi} ^{-1}_{ab} \pi_{k}^{b\, -1}
(V^{ab}_{ik})$ is a lift of $V^{ab}_{ik}$ in the same universal cover, where 
$\tilde{\phi} _{ab}$ is the map induced by $\phi _{ab}$ on the level of 
universal coverings (i.e. is the analytic representation of $\phi _{ab}$).
There exists therefore an element $\mu ^{ab}_{ik}$ of $\Lambda _D $ such that 
$\tilde{\phi} _{ab}^{-1}w^b_k(x) = \mu ^{ab}_{ik} + w^a_i(x)$ for all $x \in
V^{ab}_{ik}$. 

We now define  
\[ \psi ^{ab}_{ik} (x) = e_a ( \mu ^{ab}_{ik}, w^a_i(x),Z^a(s) ) 
\;  h^{ab} \left( \tilde{\phi} _{ab}^{-1}(w^b_k(x)), Z^a(s) \right) ^{-1}. 
\] 
It is a straight forward calculation to verify that those make the above
diagram to commute and satisfy the required relations for automorphisms.  

We define $\cl _{\rm can} := \bigsqcup _{a,i} (V^a_i  \times {\Bbb C})/ 
{\sim } $, where $(x^a_i, z^a_i) \sim (x^b_k, z^b_k)$ iff 
\[ x^a_i = x^b_k=s \;\;\;\; \mbox{and} \;\;\;\;  \left\{  \begin{array}{ll} 
z^b_k = g^a_{ik} \, z^a_i & \mbox{if $a=b$} \\ z^b_k = \psi ^{ab}_{ik} \, 
z^a_i  & \mbox{if $a \neq b$} \end{array}   \right.  \, . 
\]
Note that those satisfy the cocycle relation for transition functions. We also
have that the restriction of $\cl _{\rm can}$ to the zero section $S$  is 
trivial. Indeed, for $x \in V^{ab} \cap S$, we have that $w^a:=w^a(x), \, w^b
:=w^b(x) \in  \Lambda _D$  and $\tilde{w}^b:=\tilde{\phi} _{ab}^{-1}(w^b(x)) 
\in \Lambda _D $.  Using the defining  properties of the factors of automorphy 
and the definition of $h^{ab}$,  we get that  $e_a ( \mu ^{ab}, w^a(x),Z^a(s))
= e_a ( \tilde{w}^b, 0 ,Z^a(s) ) \; e_a ( w^a, 0 ,Z^a(s) ) ^{-1}$  and 
$h^{ab}(\tilde{w}^b, Z^a(s))^{-1} = e_b ( w^b, 0 ,Z^b(s) ) \; e_a (\tilde{w}^b,
0 ,Z^a(s) ) ^{-1} $. We conclude that for $s \in S$, we have $\psi ^{ab}(s) = 
\frac{f^a (s)}{f^b(s)}$, where $f^a : U^a \lra {\Bbb C}^{\times}$ is the 
holomorphic function given by  $ f^a(s) = e ( w^a ,0, Z^a(s) ) ^{-1}$.   

Finally observe that $\cl \cong \cl _{\rm can}$. Indeed, they both restrict to
isomorphic line bundles on the fibers (since they are both of type $D$ and  
have the same characteristic)  and they are both trivial along the zero 
section. The see-saw principle implies that they are isomorphic.

\subsection{Proof of Theorems A and B.}
\label{Proof of Theorem A.}
It suffices to prove the theorem for the line bundle $\cl := \cl _{\rm can}$
constructed in the previous section. The functions $ \theta \left[
\begin{array}{cc} c_a^1+D^{-1}n \\ c_a^2 \end{array} \right]  (v^a, Z^a(s)), \;
n \in {\Bbb Z}_D $, are theta functions for the classical factor of automorphy 
$e_a : \Lambda _D  \times {\Bbb C}^g \times \siegel \lra {\Bbb C}^{\times }$ 
of characteristic $c_a$. The line bundle $T^*_{c_a} \cl _{\frak X}|_{\pi _a
(X^a)}$ corresponds by construction to the classical factor of automorphy 
$e_a$. Let $\tilde{\sigma}_n^a$ denote the section of $T^*_{c_a} \cl _{\frak X}
|_{\pi _a(X^a)}$ corresponding to  the theta function  $ \theta \left[
\begin{array}{cc} c^1_a + D^{-1}n \\ c^2_a \end{array} \right]  (v^a, Z^a(s))$. 
Then $\sigma _n^a := \pi _a ^* \tilde{\sigma}_n^a $ is the section of  $\cl _a$
corresponding to the same theta function by identifying  $(\pi _1(X^{a}), \, 
\tilde{X}^{a}) \stackrel{p_a}{\cong} (\Lambda _D, \,{\Bbb C}^g \times  p_a
(U^{a}))$.

In the following, fix the identification $(\pi _1(X^{ab}), \, \tilde{X}^{ab}) 
\stackrel{p_a}{\cong} (\Lambda _D, \,{\Bbb C}^g \times  p_a(U^{ab}))$.
The section $\phi _{ab}^* (\tilde{\sigma}_m^b|_{\pi _b(X^{ab})})$ is a section 
of the line bundle $\phi_{ab}^*(T^*_{c_b}\cl _{\frak X} |_{\pi _b(X^{ab})})$ 
corresponding to the theta function $\phi _{ab}^* \, \theta \left[
\begin{array}{cc} c^1_b + D^{-1}m \\ c^2_b \end{array} \right]  (v^b,Z^b(s))$ 
w.r.t. the factor of automorphy $\phi_{ab}^*e_b$. We therefore get that 
the section $\pi _a ^* (\phi _{ab}^* \tilde{\sigma}_m^b |_{\pi _b(X^{ab})}) =
\pi _b^* (\tilde{\sigma}_m^b|_{\pi _b(X^{ab})})= \sigma _m^b|_{X^{ab}}$ is a 
section of $\cl _b |_{X^{ab}}$ corresponding to the same theta function. 
Let $\Phi _{ab} : \cl _a|_{X^{ab}} \lra  \cl _b|_{X^{ab}} $ be the canonical 
isomorphism defined in Section
\ref{Canonical form of line bundles on abelian fibrations.}. Then $\Phi_{ab}^*
\sigma ^b_m(s)$ corresponds to the theta function $\, h^{ab}(v^a,Z^a(s)) \, \; 
\phi _{ab}^*  \theta \left[ \begin{array}{cc} c^1_b + D^{-1}m \\ c^2_b 
\end{array} \right] (v^b, Z^b(s))$, see Section 
\ref{Factors of automorphy and line bundles.}

The set ${\cal B} _a:=<\sigma ^a_n, n \in {\Bbb Z}_D  >$ is a set of sections 
of $\cl _a$ which restricts to a base of sections of $H^0(X_s,\cl _a|_{X_s})$ 
for every  $s \in U^{ab}$, see Section 
\ref{Canonical - Classical theta functions.}. The set ${\cal B}^{'}_b:= 
< \Phi _{ab}^* \sigma ^b_m, m \in {\Bbb Z}_D > $ is also a set  of sections 
of $\cl _a|_{X^{ab}}$ which restricts to a base of sections of $H^0(X_s,\cl _a
|_{X_s})$ for every  $s \in U^{ab}$. There exists then a matrix $\tilde{C}
^{ab}(s)$ of complex functions $U^{ab} \lra {\Bbb C}^{\times }$ such that
${\cal B}^{'}_b=\tilde{C}^{ab}(s) \; {\cal B}_a$. Expressing this relation 
in terms  of theta functions, it becomes
\begin{eqnarray*}
\lefteqn{h^{ab}(v^a,Z^a(s)) <\phi _{ab}^*  \theta \left[ \begin{array}{cc} 
c^1_b + D^{-1}m  \\ c^2_b \end{array} \right]  (v^b, Z^b(s)), \; 
m \in {\Bbb Z}_D > } \\  
& &   =\tilde{C}^{ab}(s) \; <  \theta \left[ \begin{array}{cc} c^1_a + D^{-1}n 
\\ c^2_a \end{array} \right]  (v^a,Z^a(s)), \; n\in {\Bbb Z}_D >. 
\end{eqnarray*}                                                
By using  Proposition \ref{proposition1.1} and relation (\ref{eq:ttf-1}) in 
the proof of that proposition,  we conclude that $\tilde{C}^{ab}(s) =
e(\pi i \, ^tM[c]^1M[c]^2- \pi i \, ^t c^1c^2) (\det (\gamma ^{ab} Z^a(s) + 
\delta ^{ab})) ^{\frac{1}{2}}  \, C_{M^{ab}}^{-1}$  and so, 
$\det \tilde{C}^{ab}(s)\, = \zeta _8 (\det (\gamma ^{ab} Z^a(s)+\delta ^{ab})) 
^{\frac{d}{2}}$, except when $3|d_g$ and $gcd(3,d_{g-1})=1$, in which case
we have that $\det \tilde{C}^{ab}(s)\, = \zeta _{24} (\det (\gamma ^{ab} Z^a(s) 
+ \delta ^{ab})) ^{\frac{d}{2}}$.
 
By the construction of the line bundle $\cl := \cl _{\rm can}$ we have 
canonical isomorphisms $i_a :\cl |_{X^{ab}} \stackrel{\simeq}{\lra} 
\cl _a|_{X^{ab}}$ which make the following diagram commutative
\noindent
\begin{equation}
\label{diag:rotbos-1}
\begin{diagram}[aaaaaa]
\node{\cl _a |_{X^{ab}} } \arrow[2]{e,t}{\Phi _{ab}} 
\node[2]{\cl _b|_{X^{ab}}  }\\    
\node{ } \node{\cl |_{X^{ab}}} \arrow{nw,r}{i_a} \arrow{ne,r}{i_b}  
\end{diagram}
\end{equation}
\noindent  
                          
The set ${\cal B}_a^{\cl }:=<i_a^*\sigma ^a _n, n \in {\Bbb Z}_D >$ is a set 
of sections of $\cl |_{X^a}$ which restricts to a base of sections of 
$H^0(X_s,\cl |_{X_s})$ for every  $s \in U^a$. Similarly, the set ${\cal B}_b
^{\cl }:=<i_b^*\sigma ^b_m, m \in {\Bbb Z}_D >$ is a set of sections of 
$\cl |_{X^b}$ which restricts to a base of sections of $H^0(X_s,\cl |_{X_s})$ 
for every  $s \in U^b$. According now to the above Diagram 
(\ref{diag:rotbos-1}), we have on $X^{ab}={\rm f}^{-1}(U^{ab})$ that 
${\cal B}_b^{\cl }= \tilde{C}^{ab}(s) \; {\cal B}_a^{\cl }$, where 
$\tilde{C}^{ab}(s)$ is the above defined matrix. To conclude, we have proven
 
\begin{lem}
\label{lemma2.1}
Let ${\rm f}: X \lra B$ be a fibration of abelian varieties of relative 
dimension $g$ and $\cl $ the above constructed symmetric (resp. totally 
symmetric and $g \geq 3$) line bundle on $X$. Fix the trivialization $\{U^a \}$
of $B$ defined in Section
\ref{Canonical form of line bundles on abelian fibrations.}.
Then the transition functions of the line bundle $\det {\rm f}_* \cl $ are 
given by $g^{ab}_{\cl } (s) = \zeta _8  (\det (\gamma ^{ab} Z^a(s)+ \delta 
^{ab}))^{\frac{d}{2}}$ (resp. $g^{ab}_{\cl } (s) = (\det (\gamma ^{ab} Z^a(s)+
\delta ^{ab}))^{\frac{d}{2}}$), except when $3|d_g$ and $gcd(3,d_{g-1})=1$, 
in which case we have that they are given by $g^{ab}_{\cl } (s) = \zeta _{24} 
(\det (\gamma ^{ab} Z^a(s)+ \delta ^{ab}))^{\frac{d}{2}}$ (resp. $g^{ab}_{\cl }
(s) = \zeta _3 (\det (\gamma ^{ab} Z^a(s)+ \delta ^{ab}))^{\frac{d}{2}}$). 
\end{lem}

To proceed a little further we state the following lemma

\begin{lem}
\label{lemma2.2}
Let ${\rm f} : X \lra B$ be a fibration of abelian varieties and 
$s : B \lra X$ its zero section. Let $\Omega _{X/B}$ denote the relative 
cotangent bundle. Then $\Omega _{X/B} \cong {\rm f}^*E$, where $E \cong s^*
\Omega _{X/B}$ is the vector bundle on $B$ defined by the transition matrices 
$g^{ab}_E:=(\gamma ^{ab}Z^a(s)+\delta ^{ab})^{-1}$. In particular for the 
relative dualizing sheaf of ${\rm f}$ we have that $\omega _{X/B}  \cong  
{\rm f}^* \mu $, where $\mu \cong s^*\omega _{X/B}$ is the line bundle on $B$ 
defined by the transition functions  $g^{ab}_{\mu }(s)= \det (\gamma ^{ab}Z^a
(s)+\delta ^{ab})^{-1}$.  
\end{lem}    

\noindent
{\em Proof.} Let $\{ U^a \}_{a \in A}$ be an open covering of $B$ and $<U^a, 
\lambda ^a_1(s), \ldots , \lambda ^a_{2g}(s)>$ a choice of a symplectic base 
on the fibers of $X^a:= f^{-1}(U^a)$, as in Section 
\ref{Canonical form of line bundles on abelian fibrations.}. 
Let $Z^a(s)$ be the period matrix of $X^a_s$ defined by $<\lambda ^a_1(s), 
\ldots , \lambda ^a_g(s)> = <\lambda ^a_{g+1}(s), \ldots , \lambda ^a_{2g}(s)>
Z^a(s)$. $\Lambda _D$ acts on ${\Bbb C}^g \times U^a$ by $n(v,s):=(v+j_{Z^a(s)}
(n),s)$, where $j_{Z^a(s)}(n):=Z^a(s)n^1+n^2$ (where $n=(n^1,n^2)$), see 
Section  \ref{Action of the symplectic group.}. There is a canonical 
isomorphism $\phi _a : X^a \lra ({\Bbb C}^g \times U^a)/\Lambda _D $ (fibered 
over $U^a$) defined  on the level of universal coverings by sending 
$\tilde{\phi}_a(\lambda ^a_{g+i}) = (e_i,s) \in {\Bbb C}^g \times U^a, \, i=1, 
\ldots ,g$, where $<e_1, \ldots ,e_g >$ is the standard base of ${\Bbb C}^g$. 
Let $<z_1, \ldots ,z_g>$ denote the standard coordinates of ${\Bbb C}^g$. 
Then $dz_i$ is the dual to $e_i$. Let $z^a_i:=\tilde{\phi} _a^*(z_i \times
id)$. Then $<dz^a_1, \ldots , dz^a_g>$ defines at each point of $X^a$ a base 
of sections of the fiber of $\Omega _{X/B}|_{X^a}$ and $dz^a_1 \wedge \ldots 
\wedge dz^a_g$ defines a (nowhere zero) section of $\omega _{X/B}|_{X^a}$. 
This is because $dz_1 \wedge \ldots \wedge dz_g$ defines a (nowhere zero) 
section of the relative dualizing sheaf of the fibration $({\Bbb C}^g \times 
U^a)/ \Lambda _D \lra U^a$.

We have that $<\lambda ^b_1, \ldots ,\lambda ^b_{2g}> \, = \, <\lambda ^a_1, 
\ldots ,\lambda ^a_{2g}> \, ^tM^{ab}$. Therefore  
\begin{eqnarray*}
\lefteqn{<\lambda ^b_{g+1}, \ldots ,\lambda ^b_{2g}> \, = \, <\lambda ^a_1, 
\ldots ,\lambda ^a_{2g}> \left(  \begin{array}{c} \, ^t \gamma ^{ab}\\ \, ^t 
\delta ^{ab} \end{array} \right)= } \\ 
& & <\lambda ^a_{g+1}, \ldots ,\lambda ^a_{2g}>(Z^a(s),I) \, \left( 
\begin{array}{c} \, ^t \gamma ^{ab} \\ \, ^t \delta ^{ab} \end{array} \right) 
= \; <\lambda ^a_{g+1}, \ldots ,\lambda ^a_{2g}> \, ^t (\gamma ^{ab} Z^a(s) + 
\delta ^{ab}).
\end{eqnarray*}  
\noindent
By taking dual bases we have that 
\[ <dz ^b_1, \ldots , dz ^b_g > \, = \, <dz ^a_1, \ldots , dz ^a_g > \, 
(\gamma ^{ab} Z^a(s) + \delta ^{ab} )^{-1}
\]
and so, 
\[ dz ^b_1 \wedge \ldots \wedge  dz ^b_g \, = \, {\rm det } ( \gamma ^{ab} Z^a
(s) + \delta ^{ab} ) ^{-1} \, dz ^a_1 \wedge \ldots \wedge  dz ^a_g. 
\]  
This proves the lemma.
\begin{flushright}
${\Box}$
\end{flushright} 

The proof of Theorems \ref{theorem1} and \ref{theorem2} is a consequence 
of Lemma \ref{lemma2.1} and Lemma \ref{lemma2.2}.

\subsection{The Jacobian fibration - Proof of Theorem C.}
\label{Proof of Theorem C.} 
We apply now the above considerations to the Jacobian fibration ${\rm f} : 
\jac  \lra \mod $, where $\jac $ denotes the universal Jacobian
parametrizing line bundles of degree zero on the fibers of the universal
curve $\psi : {\cal C} \lra \mod $ over the moduli space $\mod $ 
of smooth, irreducible curves of genus $g \geq 3$ without automorphisms. 
The Picard group of $\mod $ is generated by a line bundle $\lambda $ which is
defined to be the determinant of the Hodge bundle of the map $\psi $ i.e.
$\lambda := \det \psi_* \omega_{{\cal C}/ \mod  }$, see \cite{AC}.

On the Jacobian fibration $\jac $, there is a totally symmetric line bundle 
$\cl $ which restricts to a line bundle of class $2 \theta $ on the fibers 
and it is trivial along the zero section. It is defined as the pull back of the 
Poincare bundle under the natural map. By  Theorem  \ref{theorem2}, we get that

\begin{cor}
\label{corollary2.1}  
Let $\cl $ denote the above canonical line bundle on $\jac $. Then
\[ \det {\rm f}_*(\cl ^{\otimes n}) \cong - \frac{(2n)^g}{2} 
\det \psi _*\omega_{{\cal C}/ \mod }. 
\]
\end{cor}

\noindent
\begin{rem} {\rm  One can also prove the above Corollary \ref{corollary2.1}
by using the Theorem 5.1 in \cite{FC}. For the case of the universal Jacobian,
we can ged rid from the torsion factor which appears in that theorem, by using
the fact that the Picard group of $\mod $ is generated by 
$\det \psi _*\omega_{{\cal C}/ \mod }$, see \cite{AC}. }
\end{rem} 

Let ${\rm f}_{g-1} : \jac ^{g-1} \lra \mod $ be the Jacobian fibration of 
degree $g-1$. We use the following result from \cite{T1}, \cite{T2}.
Let $\alpha: \tilde{\mod } \lra \mod$ denote the covering of even theta
characteristics in $\jac ^{g-1}$. It is a covering of degree $2^{g-1}(2^g+1)$.
The theta divisor $\Theta $ in $\jac ^{g-1}$ intersects $\tilde{\mod }$ 
transversely and the (set theoretic) intersection projects birationally via 
$\alpha $ to a divisor in  $\mod $ which has class  $2^{g-3}(2^g+1)\,{\rm c}_1
( \lambda )$.  On the other hand the generic point of the intersection 
corresponds to a line bundle with two sections. By the description of the 
singularities of the theta divisor, we have that on such a point the theta 
divisor has a singularity of  multiplicity $2$. We therefore have that the 
push-forward via $\alpha $ of the (scheme theoretic) intersection of $\Theta $ 
with $\tilde{\mod }$ is a divisor of class $2^{g-2}(2^g+1)\, {\rm c}_1
( \lambda )$.

We use the following commutative diagram. 
\noindent
\begin{equation}
\label{diagram2.5}
\begin{diagram}[aaaa]
\node{\jac } \arrow{s,l}{{\rm f}} 
\node{\tilde{\jac }} \arrow{w,t}{\gamma } \arrow{se,r}{\tilde{{\rm f}}} 
\arrow[2]{e,t}{\phi }  
\node[2]{\tilde{\jac }^{g-1}}  \arrow{sw,r}
{\tilde{{\rm f}}_{g-1}}  \arrow{e,t}{\delta }
\node{\jac ^{g-1}} \arrow{s,r}{{\rm f}_{g-1}}  \\
\node{\mod }  
\node[2]{\tilde{\mod }} \arrow[2]{w,t}{\alpha } \arrow[2]{e,t}{\alpha }
\node[2]{\mod } 
\end{diagram}
\end{equation}

\noindent   
In the above diagram we denote by $\tilde{\jac }$ and $\tilde{\jac } ^{g-1}$ 
the pull back of $\jac $ and $\jac ^{g-1}$ on $\tilde{\mod }$. By $\phi $ we
denote the map which sends the line bundle $L \in \tilde{\jac }$ sitting 
over the point $[C, \eta ] \in \tilde{\mod }$ to the line bundle  
$L ^{\otimes 2} \otimes \eta  \in \tilde{\jac }^{g-1}$. The map $\phi $ is 
etale of degree $2^{2g}$. Let $s: \tilde{\mod } \lra \tilde{\jac }$ denote the 
zero section and let $\sigma :  \tilde{\mod } \lra \tilde{\jac }^{g-1}$
denote the section which sends $[C, \eta ] \ms \eta$. We have that $\phi \,
s = \sigma $. 

Let $\tilde{\Theta } $ be the line bundle corresponding to the theta divisor 
on $\tilde{\jac }^{g-1}$. Then $\tilde{\Theta } =  \delta ^* \Theta $ and so, 
$\alpha _* {\rm c}_1( \tilde{\rm f}_{g-1 \; *} \tilde{\Theta } ^{\otimes n} ) 
= 2^{g-1}(2^g+1) {\rm c}_1( {\rm f}_{g-1 \; *} \Theta ^{\otimes n} ) $. Let 
$\tilde{\lambda }$ be the determinant of the Hodge bundle of the fibration 
$\tilde{\rm f}$. Then $\tilde{\lambda} = \alpha ^* \lambda$ and so $\alpha _* 
{\rm c}_1(\tilde{\lambda }) = 2^{g-1}(2^g+1) {\rm c}_1({\lambda })$. 
From the things we stated above, we have that if $\tilde{\mu} :=\sigma ^* 
\tilde{\Theta }$ then $\alpha _* {\rm c}_1(\tilde{\mu }) = 2^{g-2}(2^g+1)\, 
{\rm c}_1( \lambda )$. Let $\tilde{\cl }$ be the canonical line bundle 
on $\tilde{\jac }$ of Corollary \ref{corollary2.1}. Then $\tilde{\cl } = 
\gamma ^* \cl $. 

One can see that ${\rm c}_1(\tilde{\rm f}_* \phi ^* \tilde{\Theta }^{\otimes 
n}) = 2^{2g} {\rm }c_1 (\tilde{\rm f}_{g-1 \; *} \tilde {\Theta }^{\otimes n}
)$. A proof can be given by applying the GRR theorem. One can also see that the
restriction of $\phi ^*\tilde{\Theta }$ on a fiber of the map $\tilde{{\rm f}}$
is equal to the restriction of $\tilde{\cl }^{\otimes 2}$ on the same fiber. 
This can be seen by using Proposition 3.5 of ${\rm Ch.} \, 2$ in \cite{LB} and 
the Riemann's constant theorem. Therefore, by the see-saw principle, 
the line bundles $\tilde{\cl }^{\otimes 2}$ and $\phi ^* \tilde{\Theta }$  
are isomorphic up to tensor by a line bundle which is pull back of a line 
bundle from $\tilde{\mod }$. Since $s^* \tilde{\cl }^{\otimes 2} \cong 
{\cal O}$ and $s^* \phi ^* \tilde{\Theta } \cong \tilde{\mu }$, we have that
\[ \tilde{\cl }^{\otimes 2n} \; \otimes \; \tilde{{\rm f}}^* \tilde{\mu }
^{\otimes n} \; \cong \; \phi ^*\tilde{\Theta }^{\otimes n}.      
\] 
By applying $\tilde{{\rm f}}_* $ and taking ${\rm c}_1$ we have that 
\[ {\rm c}_1( \tilde{{\rm f}}_*\tilde{\cl }^{\otimes 2n} ) \, + \,
(4n)^gn {\rm c}_1( \tilde{\mu } ) \, = \, 2^{2g}{\rm c}_1
( \tilde{{\rm f}}_{g-1 \; *} \tilde{\Theta }^{\otimes n} ).     
\]
Apply now $\alpha _*$ to get
\[ -2^{g-1}(2n)^g \,2^{g-1}(2^g+1) {\rm c}_1( \lambda ) + (4n)^g n \, 
2^{g-2}(2^g+1) {\rm c}_1( \lambda ) \, = \, 2^{2g}\, 2^{g-1}(2^g+1) 
{\rm c}_1 ( {\rm f}_{g-1 \; *} \Theta ^{\otimes n} ).
\]
Therefore
\[ {\rm c}_1({\rm f}_{g-1 \; *} \Theta ^{\otimes n }) =
\frac{1}{2} n^g (n-1) {\rm c}_1(\lambda ).
\]
By using the fact that the ${\rm Pic}\mod $ is generated by $\lambda $, see
\cite{AC}, this concludes the proof of the Theorem \ref{theorem2}.

\subsection{Alternative proof of Theorem C.}
\label{Alternative proof of Theorem C.}
We give now the alternative proof of Theorem \ref{theorem3} promised in the
introduction. This is an application of the GRR theorem, see also Appendice 2
in \cite{MB1} for a similar calculation. We keep the notation of the above 
section \ref{Proof of Theorem C.}. In the above Diagram  \ref{diagram2.5} let 
$\phi $ be the map which sends the line bundle $L \in \tilde{\jac }$ sitting 
over the point $[C, \eta ] \in \tilde{\mod }$ to the line bundle  $L \otimes 
\eta  \in \tilde{\jac }^{g-1}$. By Lemma  \ref{lemma2.2} we have that 
$\Omega _{\tilde{\jac } / \tilde{\mod}} \cong \tilde{\rm f} ^* \tilde{s}^* 
\Omega _{\tilde{\jac } / \tilde{\mod}}$ and since $\phi $ is an isomorphism, 
we get that $\Omega _{\tilde{\jac }^{g-1} / \tilde{\mod}} \cong \tilde{\rm f}
_{g-1} ^* \tilde{s}^* \Omega _{\tilde{\jac }^{g-1} / \tilde{\mod}}$. 
     
We proceed now to the proof of Theorem \ref{theorem3}: We apply GRR to the 
fibration $\tilde{\rm f}_{g-1} :\tilde{\jac } ^{g-1} \lra  \tilde{\mod }$.
It is
\[ {\rm ch}(\tilde{\rm f}_{g-1 \; !} (\tilde{\Theta } ^{\otimes n})) =  
\tilde{\rm f}_{g-1 \; *}({\rm ch} (\tilde{\Theta }^{\otimes n}) \cdot 
{\rm td} (\Omega ^{\vee} _{\tilde{\jac }^{g-1} / \tilde{\mod}})).
\]
We get that 
\[ {\rm c_1}(\tilde{f}_{g-1 \; *} \tilde{\Theta }^{\otimes n}) 
= \frac{n^{g+1}}{(g+1)!} f_*{\rm c_1}^{g+1}(\tilde{\Theta }) - 
\frac{n^g}{2}{\rm c_1}( \tilde{\lambda }). 
\]
The vanishing of the terms on the right hand side containing the ``factor"
${\rm c_1}^k$ with $k \leq g-1$ in the expansion of ${\rm ch}(\tilde{\Theta }
^{\otimes n})$, is a consequence of the projection formula and the fact that 
$\Omega _{\tilde{\jac }^{g-1} / \tilde{\mod}} \cong \tilde{\rm f}_{g-1} ^* E $,
where $E$ is a vector bundle, see Lemma \ref{lemma2.2}. The form of the 
term containing the ``factor" ${\rm c_1}^{g-1}$ is due to the Poincare 
formula. The appearance of $\tilde{\lambda}$ is a consequence of 
Corollary  \ref{corollary2.1}.
                             
Lets say now that ${\rm c_1}({\rm f}_{g-1 \; *} \Theta ^{\otimes n}) 
= a(n){\rm c_1} ( \lambda )$ and  that ${\rm f}_*{\rm c_1} ^{g+1}(\Theta ) = 
b {\rm c_1}(\lambda )$, where $a(n) , b \in {\Bbb Z}$, see  \cite{AC}. Then 
${\rm c_1}(\tilde{\rm f}_{g-1 \; *} \tilde{\Theta } ^{\otimes n}) = 
a(n) {\rm c_1}( \tilde{\lambda })$  and  $\tilde{\rm f}_*{\rm c_1} ^{g+1}
(\Theta ) =  b {\rm c_1}(\tilde{\lambda })$. We get  $a(n) = 
\frac{n^{g+1}}{(g+1)!} b - \frac{n^g}{2}$. For $n=1$ the above gives that 
$b= (g+1)! (a(1)+\frac{1}{2})$. Now $a(1)=0$, because the line bundle 
$\tilde{\rm f}_{g-1 \; *} \tilde{\Theta }$ has by definition a nowhere zero 
section, so it is the trivial bunlde. We then get  $b=\frac{(g+1)!}{2}$ and so, 
${\rm c_1}({\rm f}_{g-1 \; *} \Theta ^{\otimes n}) = \frac{1}{2}n^g(n-1)
{\rm c_1}(\lambda )$.

\end{document}